\documentclass{aims}
\usepackage{amsmath}
  \usepackage{paralist}
  \usepackage{graphics} 
  \usepackage{epsfig} 
\usepackage{graphicx}  \usepackage{epstopdf}
 \usepackage[colorlinks=true]{hyperref}
\hypersetup{urlcolor=blue, citecolor=red}

\def\currentvolume{X}
 \def\currentissue{X}
  \def\currentyear{200X}
   \def\currentmonth{XX}
    \def\ppages{X--XX}

\def\R{\mathbb R}

\def\N{\mathbb N}

\def\romega{{r_{\phantom{a}\!\!\!\!_\Omega}}}
\def\rgamma{{r_{\phantom{a}\!\!\!\!_\Gamma}}}
\def\loc{{\text{\upshape loc}}}
\newcommand{\cal}[1]{{\mathcal #1}}
\DeclareMathOperator\essinf{essinf}

\newtheorem{theorem}{Theorem}[section]

\newtheorem{lemma}[theorem]{Lemma}

\theoremstyle{definition}
\newtheorem{definition}[theorem]{Definition}
\newtheorem{remark}{Remark}

\title[Blow--up for the wave equation with hyperbolic...]
{Blow--up for the  wave equation with hyperbolic dynamical  boundary
conditions, interior and boundary nonlinear damping and sources}

\author[Enzo Vitillaro]{}


 \email{enzo.vitillaro@unipg.it}

\thanks{The work was realized within the auspices of the INdAM -- GNAMPA Projects
{\em Equazioni alle derivate parziali: Problemi e Modelli} (Prot\_U-UFMBAZ-2020-000761), and it was also supported by {\em Progetto Equazione delle onde con condizioni acustiche,  finanziato  con  il Fondo  Ricerca  di Base, 2019, della Universit\`a degli Studi di Perugia} and by {\em Progetti Equazioni delle onde con condizioni iperboliche ed acustiche al bordo,  finanziati  con  i Fondi  Ricerca  di Base 2017 and 2018, della Universit\`a degli Studi di Perugia}.}

\begin{document}
\maketitle

\centerline{\scshape Enzo Vitillaro}
\medskip
{\footnotesize
 \centerline{Dipartimento di Matematica ed Informatica, Universit\`a di Perugia}
      \centerline{Via Vanvitelli,1 06123 Perugia ITALY}
} 

\begin{abstract}
The aim of this paper is to give global nonexistence and blow--up results for  the problem
$$
\begin{cases} u_{tt}-\Delta u+P(x,u_t)=f(x,u) \qquad &\text{in
$(0,\infty)\times\Omega$,}\\
u=0 &\text{on $(0,\infty)\times \Gamma_0$,}\\
u_{tt}+\partial_\nu u-\Delta_\Gamma u+Q(x,u_t)=g(x,u)\qquad
&\text{on
$(0,\infty)\times \Gamma_1$,}\\
u(0,x)=u_0(x),\quad u_t(0,x)=u_1(x) &
 \text{in $\overline{\Omega}$,}
\end{cases}$$
where $\Omega$ is a bounded  open $C^1$ subset of $\R^N$, $N\ge 2$,
$\Gamma=\partial\Omega$, $(\Gamma_0,\Gamma_1)$ is a  partition of $\Gamma$,  $\Gamma_1\not=\emptyset$ being relatively open in $\Gamma$, $\Delta_\Gamma$ denotes the
Laplace--Beltrami operator on $\Gamma$, $\nu$ is the outward normal
to $\Omega$, and the terms $P$ and $Q$ represent nonlinear damping
terms, while $f$ and $g$ are nonlinear source terms.
These results complement the analysis of the problem given by the author in two recent papers, dealing with local and global existence, uniqueness and well--posedness.
\end{abstract}

\section{Introduction and main results} \label{intro}
\subsection{Presentation of the problem and literature overview}

We deal with the evolution problem consisting of
the
 wave equation posed in a bounded regular open subset  of
$\R^N$, supplied with a second order dynamical boundary condition of
hyperbolic type, in presence of interior and/or boundary damping
terms and sources. More precisely we consider the
initial --and--boundary value problem
\begin{equation}\label{1.1}
\begin{cases} u_{tt}-\Delta u+P(x,u_t)=f(x,u) \qquad &\text{in
$(0,\infty)\times\Omega$,}\\
u=0 &\text{on $(0,\infty)\times \Gamma_0$,}\\
u_{tt}+\partial_\nu u-\Delta_\Gamma u+Q(x,u_t)=g(x,u)\qquad
&\text{on
$(0,\infty)\times \Gamma_1$,}\\
u(0,x)=u_0(x),\quad u_t(0,x)=u_1(x) &
 \text{in $\overline{\Omega}$,}
\end{cases}
\end{equation}
where $\Omega$ is a bounded $C^1$ open subset of $\R^N$, with $N\ge
2$. We denote  $\Gamma=\partial\Omega$ and we assume that
$\Gamma=\Gamma_0\cup\Gamma_1$, $\Gamma_0\cap\Gamma_1=\emptyset$,
that $\Gamma_1\not=\emptyset$ is relatively open in $\Gamma$ and, denoting by
$\sigma$ the  standard Lebesgue hypersurface measure on $\Gamma$,  that
$\sigma(\overline{\Gamma}_0\cap\overline{\Gamma}_1)=0$.
These properties of $\Omega$, $\Gamma_0$ and $\Gamma_1$ will be assumed,
without further comments, throughout the paper.  Moreover $u=u(t,x)$, $t\ge 0$, $x\in\Omega$,
$\Delta=\Delta_x$ denotes the Laplace operator respect to the space
variable, while $\Delta_\Gamma$ denotes the Laplace--Beltrami
operator on $\Gamma$ and $\nu$ is the outward normal to $\Omega$.

The terms $P$ and $Q$ represent nonlinear damping terms, i.e.
$P(x,v)v\ge 0$, $Q(x,v)v\ge 0$,  the cases $P\equiv 0$ and $Q\equiv
0$ being specifically allowed, while  $f$ and $g$  represent
nonlinear source terms. The specific assumptions on them
will be introduce later on.

Problems with kinetic boundary conditions, that is boundary
conditions involving $u_{tt}$ on $\Gamma$, or on a part of it,
naturally arise in several physical applications. A one dimensional
model was studied by several authors to describe transversal small
oscillations of an elastic rod with a tip mass on one endpoint,
while the other one is pinched. See
\cite{andrewskuttlershillor,conradmorgul,darmvanhorssen,guoxu,markuslittman1,markuslittman2,morgulraoconrad} and also
\cite{meurerkugi} were a piezoelectric stack actuator is modeled.

 A two dimensional model introduced in
\cite{goldsteingiselle} deals with a vibrating membrane of surface
density $\mu$, subject to a tension $T$, both taken constant and
normalized here for simplicity.  If $u(t,x)$,
$x\in\Omega\subset\R^2$  denotes the vertical displacement from the
rest state, then (after a standard linear approximation) $u$
satisfies the wave equation $u_{tt}-\Delta u=0$,
$(t,x)\in\R\times\Omega$. Now suppose that a part $\Gamma_0$ of the
boundary  is pinched, while the other part $\Gamma_1$ carries a
constant linear mass density $m>0$ and  it is subject to a linear
tension $\tau$. A practical example of this situation is given by a
drumhead with a hole in the interior having  a thick border, as
common in bass drums. One linearly approximates the force exerted by
the membrane on the boundary with $-\partial_\nu u$. The boundary
condition thus reads as $mu_{tt}+\partial_\nu u-\tau
\Delta_{\Gamma_1}u=0$. In the quoted paper the case
$\tau=0$ was studied (when $\Gamma_0=\emptyset$), while here we
consider the more realistic case $\tau>0$, with $\tau$ and $m$ normalized for simplicity, and we also allow
$\Gamma_0$ to be nonempty. We would
like to mention that this model belongs to a more general class of
models of Lagrangian type involving boundary energies, as introduced
for example in \cite{fagotinreyes}.

A three dimensional model involving kinetic dynamical boundary
conditions comes out from \cite{GGG}, where a gas undergoing small
irrotational perturbations from rest in a domain $\Omega\subset\R^3$
is considered. Normalizing the constant speed of propagation, the
velocity potential $\phi$ of the gas (i.e. $-\nabla \phi$ is the
particle velocity) satisfies the wave equation $\phi_{tt}-\Delta
\phi=0$ in $\R\times\Omega$. Each point $x\in\partial\Omega$ is
assumed to react to the excess pressure of the acoustic wave like a
resistive harmonic oscillator or spring, that is the boundary is
assumed to be locally reacting (see \cite[pp.
259--264]{morseingard}). The normal displacement $\delta$ of the
boundary into the domain then satisfies
$m\delta_{tt}+d\delta_t+k\delta+\rho\phi_t=0$, where $\rho>0$ is the
fluid density  and $m,d,k\in C(\partial\Omega)$, $m,k>0$, $d\ge 0$.
When the boundary is nonporous one has $\delta_t=\partial_{\nu}\phi$
on $\R\times\partial\Omega$, so the boundary condition reads as
$m\delta_{tt}+d\partial_\nu \phi+k\delta+\rho\phi_t=0$. In the
particular case $m=k$ and $d=\rho$  (see  \cite[Theorem~2]{GGG}) one
proves that $\phi_{|\Gamma}=\delta$, so the boundary condition reads
as  $m\phi_{tt}+d\partial_\nu \phi+k\phi+\rho\phi_t=0$, on
$\R\times\partial\Omega$.
 Now, if one consider the case in which the
boundary is not locally reacting, as in \cite{beale}, one adds
a Laplace--Beltrami term  so getting a dynamical
boundary condition like in \eqref{1.1}. See \cite{mugnvit} where this case was studied in detail.

Several papers in the literature deal with the wave equation  with
kinetic boundary conditions. This fact is even more evident  if one
takes into account that, plugging the equation in \eqref{1.1} into
the boundary condition, we can rewrite it as  $\Delta u
+\partial_\nu u-\Delta_\Gamma u+ Q(x,u_t)+P(x,u_t)=f(x,u)+g(x,u)$.
Such a condition is usually called a {\em generalized Wentzell
boundary condition}, at least when nonlinear perturbations are not
present. We refer to
\cite{doroninlarkinsouza,FGGGR,littmanliu,mugnolo,vazvitM3AS,xiaoliang2,xiaoliang1,zhang}.
All of them deal  either with the case $\tau=0$ or with linear
problems.

Here we shall consider this type of kinetic boundary condition in
connection with nonlinear boundary damping and source terms. These
terms have been considered by several authors, but mainly in
connection with first order dynamical boundary conditions. See
\cite{MR2674175,
MR2645989,bociuNAMA,bociulasiecka2,bociulasiecka1,CDCL,CDCM,chueshovellerlasiecka,fisvit2,lastat,global, zuazua}.
The competition between interior damping and source terms is
methodologically related to the competition between boundary damping
and source and it possesses a large literature as well. See
\cite{MR2609953,georgiev,levserr,ps:private,radu1,STV,blowup}.

Local and global existence, continuation, uniqueness and Hadamard well--posedness for problem \eqref{1.1} has been studied by  the author in the
recent papers \cite{Dresda1, Dresda2} (see also \cite{AMS} for a preliminary study of a particular case). In \cite{Dresda1} a blow--up result was also given when $P$ and $Q$ are linear in $u_t$.

Moreover a linear problem strongly related to \eqref{1.1} has also
been recently studied in \cite{Fourrier, graberlasiecka}, and another one in the recent paper \cite{zahn}, dealing with
holography, a main theme in theoretical high energy physics and
quantum gravity. See also \cite{graber}.

The aim of the present paper is to discuss the optimality of the global existence result in \cite{Dresda2} by giving some complementary global nonexistence and blow--up results for solutions of \eqref{1.1} when $P$ and $Q$ are possible nonlinear in $u_t$, a case which remained open in \cite{Dresda1}.

To simplify the presentation of our main results we shall restrict, in this section, to a parameters--dependent family of model problems, which catches the essential features of \eqref{1.1}, as long as the alternative between global existence and nonexistence for arbitrary initial data is concerned.
\subsection{A family of model problems} We shall deal with
\begin{equation}\label{1.2}
\begin{cases} u_{tt}-\Delta u+\alpha \left(a|u_t|^{\widetilde{m}-2}u_t+\!\!|u_t|^{m-2}u_t\right)=\gamma |u|^{p-2}u  &\text{in
$(0,\infty)\times\Omega$,}\\
u=0 &\text{on $(0,\infty)\times \Gamma_0$,}\\
u_{tt}+\partial_\nu u-\!\Delta_\Gamma u+\!\beta \left(b|u_t|^{\widetilde{\mu}-2}u_t+\!\!|u_t|^{\mu-2}u_t\right)=\delta |u|^{q-2}u
&\text{on
$(0,\infty)\times \Gamma_1$,}\\
u(0,x)=u_0(x),\quad u_t(0,x)=u_1(x) &
 \text{in $\overline{\Omega}$,}
\end{cases}
\end{equation}
where $a,b,\alpha,\beta,\gamma,\delta,\widetilde{m},m,\widetilde{\mu},\mu, p,q$ are real number verifying
\begin{equation}\label{1.3}
a,b,\alpha,\beta,\gamma,\delta\ge 0,\qquad 1<\widetilde{m}\le m,\qquad 1<\widetilde{\mu}\le \mu,\qquad p,q\ge 2.
\end{equation}
The terms $a|u_t|^{\widetilde{m}-2}u_t$ and $b|u_t|^{\widetilde{\mu}-2}u_t$ are present only for modeling purpose and they need a suitable handling, but their possible vanishing is not relevant in the subsequent discussion, so the reader can take $a=b=0$ in the sequel.

By the contrary the possible vanishing of each parameter among $\alpha$, $\beta$, $\gamma$ and $\delta$ individuates a different model problem in the family, which is then constituted by sixteen (!)  different model problems. The unitary treatment of them was a characteristic feature of \cite{Dresda1, Dresda2} but, when dealing with the alternative between global existence and nonexistence it is useful to introduce a classification.
In doing it we shall use the standard terminology widely adopted in the literature when dealing with a strongly methodologically related family of model problems, i.e. (taking $a=b=0$ for simplicity)
\begin{equation}\label{1.4}
\begin{cases} u_{tt}-\Delta u+\alpha |u_t|^{m-2}u_t=\gamma |u|^{p-2}u \qquad &\text{in
$(0,\infty)\times\Omega$,}\\
u=0 &\text{on $(0,\infty)\times \Gamma_0$,}\\
\partial_\nu u+u+\beta |u_t|^{\mu-2}u_t=\delta |u|^{q-2}u\qquad
&\text{on
$(0,\infty)\times \Gamma_1$,}\\
u(0,x)=u_0(x),\quad u_t(0,x)=u_1(x) &
 \text{in $\overline{\Omega}$,}
\end{cases}
\end{equation}
where  $\alpha,\beta,\gamma,\delta,\widetilde{m},m,\widetilde{\mu},\mu, p,q$ are as before.
We shall often refer to literature concerning \eqref{1.4}.

Our \emph{primary classification} concerns the presence of the interior source $|u|^{p-2}u$ and of the boundary source $|u|^{q-2}u$ (the constants $\alpha,\beta,\gamma,\delta$ could be normalized when positive), so defining the following four classes of model problems:
\renewcommand{\labelenumi}{\Alph{enumi})}
\begin{enumerate}
  \item \emph{sourceless}, when $\gamma=\delta=0$;
  \item \emph{with boundary source}, when $\gamma=0<\delta$;
  \item \emph{with interior source}, when $\delta=0<\gamma$;
  \item \emph{with interior and boundary sources}, when $\gamma,\delta>0$.
\end{enumerate}
Clearly each class includes four different model problems according to the possible presence of interior and boundary damping terms. We shall set this one to be our \emph{secondary classification}:
\renewcommand{\labelenumi}{\alph{enumi})}
\begin{enumerate}
  \item \emph{undamped}, when $\alpha=\beta=0$;
  \item \emph{with boundary damping}, when $\alpha=0<\beta$;
  \item \emph{with interior damping}, when $\beta=0<\alpha$;
  \item \emph{with interior and boundary damping}, when $\alpha,\beta>0$.
\end{enumerate}
When referring to a specific model we shall sometimes use the two letters classification given by previous lists. For example Cb) stands for the problem with interior source and boundary damping.

Source terms are usually classified in the literature concerning \eqref{1.4} (see \cite{bociulasiecka2, bociulasiecka1}), and also in \cite{Dresda2}, according to the relation
occurring between their growth and the critical exponents $\romega$ and $\rgamma$ of the Sobolev Embeddings of $H^1(\Omega)$ and
$H^1(\Gamma)$ into the corresponding Lebesgue spaces, i.e.
\begin{equation}\label{1.5}
\romega=
\begin{cases}
\dfrac {2N}{N-2} &\text{if $N \ge 3$,}\\ \infty &\text{if $N=2$},
\end{cases}
\qquad \rgamma=
\begin{cases}
\dfrac {2(N-1)}{N-3} &\text{if $N \ge 4$,}\\ \infty &\text{if $N=2,
3$}.
\end{cases}
\end{equation}
In particular the source $\gamma |u|^{p-2}u$ is:
\renewcommand{\labelenumi}{{(\roman{enumi})}}
\begin{enumerate}
\item {\em  subcritical} if $\gamma=0$ or $2\le p\le 1+\romega/2$, when the
Nemitskii operator $\widehat{\gamma |u|^{p-2}u}$ is locally
Lipschitz from $H^1(\Omega)$ into $L^2(\Omega)$;
\item {\em supercritical} if $\gamma>0$ and $1+\romega/2<p\le \romega$, when  $\widehat{\gamma |u|^{p-2}u}$ is no longer locally Lipschitz
from $H^1(\Omega)$ into $L^2(\Omega)$ but it still possesses a potential energy in
$H^1(\Omega)$;
\item {\em super--supercritical} if $\gamma>0$ and $p>\romega$, when $\widehat{\gamma |u|^{p-2}u}$ has no potentials in $H^1(\Omega)$.
\end{enumerate}
The  analogous classification is made for $\delta |u|^{q-2}u$ depending on $\delta$, $q$ and $\rgamma$.

In \cite{Dresda1} we studied well--posedness of \eqref{1.2} when both sources are subcritical, while in \cite{Dresda2} this condition was relaxed. In the sequel we shall deal with weak solutions of \eqref{1.2}, already introduced in the quoted papers. They are solutions in a suitable distribution sense and enjoy "good properties", see Definition~\ref{definition2.2} and Lemma~\ref{lemma2.3}. An essential ingredient in their definition is that
\begin{footnote}
{The sets in the planes $(p,m)$ and $(q,\mu)$, for which \eqref{1.6}
holds, corresponding to the classification above, are illustrated in
dimensions $N=2,3,4$ in Figure~1. Clearly the two sets are respectively relevant only when $\gamma>0$ and $\delta>0$, and\eqref{1.6} can be disregarded when $N=2$.}
\end{footnote}
 \begin{equation}\label{1.6}
 p\le
\begin{cases}
1+\romega/2&\text{if $\gamma>0$, $\alpha=0$,}\\
1+\romega/ \overline{m}'&\text{if $\gamma>0$, $\alpha>0$,}
\end{cases}\quad
q\le
\begin{cases}
1+\rgamma/2&\text{if $\delta>0$, $\beta=0$,}\\
1+\rgamma/ \overline{\mu}'&\text{if $\delta>0$, $\beta>0$,}
\end{cases}
\end{equation}
where, for any $\rho\in [1,\infty]$ we denote by $\rho'$ its H\"{o}lder conjugate of $\rho$, i.e. $1/\rho+1/\rho'=1$, and
\begin{equation}\label{1.7}
\overline{m}=\max\{2,m\},\qquad \overline{\mu}=\max\{2,\mu\}.
\end{equation}
The notation \eqref{1.7} will be consistently used throughout the paper.

Solutions of \eqref{1.2} in the sense of distributions may be considered also when \eqref{1.6} does not hold. Beside the lack of an available local existence theory,
any discussion on the life--span of these solutions  looks to be out of reach.
\subsection{Known results} \label{subsection1.3} In the paper we shall
identify $L^\rho(\Gamma_1)$, $\rho\in[1,\infty]$,  with its isometric image in $L^\rho(\Gamma)$,
that is
\begin{equation}\label{1.8}
L^\rho(\Gamma_1)=\{u\in L^\rho(\Gamma): u=0\,\,\text{a.e. on
}\,\,\Gamma_0\}.
\end{equation}
We shall denote by $u_{|\Gamma}$ the trace
on $\Gamma$ of $u\in H^1(\Omega)$. We introduce the Hilbert
spaces $H^0 = L^2(\Omega)\times L^2(\Gamma_1)$ and
\begin{equation}\label{1.9}
H^1 = \{(u,v)\in H^1(\Omega)\times H^1(\Gamma): v=u_{|\Gamma}, v=0
\,\,\ \text{a.e. on $\Gamma_0$}\},
\end{equation}
with the standard product norm. For the sake of
simplicity we shall identify, when useful, $H^1$ with its isomorphic
counterpart $\{u\in H^1(\Omega): u_{|\Gamma}\in H^1(\Gamma)\cap
L^2(\Gamma_1)\}$, through the identification $(u,u_{|\Gamma})\mapsto
u$, so we shall write, without further mention, $u\in H^1$ for
functions defined on $\Omega$. Moreover we shall drop the notation
$u_{|\Gamma}$, when useful, so we shall write $\|u\|_{2,\Gamma}$,
$\int_\Gamma u$, and so on, for $u\in H^1$.

We shall also use the main phase space for problem \eqref{1.2}, that is
\begin{equation}\label{1.10}
\cal{H}=H^1\times H^0,\quad\text{with the standard norm}\quad \|(u,v)\|_{\cal{H}}^2=\|u\|_{H^0}^2+\|v\|_{H^1}^2.\end{equation}
As a particular case of \cite[Theorems~1.1 and 1.2]{Dresda2} (see Theorems~\ref{theorem5.1}--\ref{theorem5.2}  below when $p\le \romega$ and $q\le \rgamma$) the following results hold. When \eqref{1.3}, \eqref{1.6} hold and
\begin{equation}\label{1.11}
\begin{alignedat}2
  &p<1+\romega/m'\qquad&&\text{when $N\ge 5$,\quad $\gamma>0$, \quad $m>\romega$},\\
  &q<1+\rgamma/\mu'\qquad&&\text{when $N\ge 6$, \quad $\delta>0$, \quad $\mu>\rgamma$,}
\end{alignedat}
\end{equation}
for all $U_0:=(u_0,u_1)\in\cal{H}$ such that
\begin{equation}\label{1.12}
\begin{alignedat}2
 &u_0\in L^{\romega(p-2)/(\romega-2)} && \quad\text{if $N=3,4$, \quad $\gamma>0$, \quad $p=1+\romega/m'$, \quad $m>\romega$, } \\
 &{u_0}_{|\Gamma}\in L^{\rgamma(q-2)/(\rgamma-2)} &&\quad \text{if $N=4,5$, \quad $\delta>0$, \quad $q=1+\rgamma/\mu'$, \quad $\mu>\rgamma$,}
 \end{alignedat}
\end{equation}
problem \eqref{1.2} possesses a maximal weak solution
$u\in C([0,T_{\text{max}});H^1)\cap C^1([0,T_{\text{max}});H^0)$ for some $T_{\text{max}}\in (0,\infty]$.
In the sequel we shall  denote $U=(u,u')\in C([0,T_{\text{max}});\cal{H})$.

It is worth observing that \eqref{1.11} and \eqref{1.12} can be disregarded when $p\le \romega$ and $q\le \rgamma$, since $m>\romega$ and $\mu>\rgamma$ respectively yield $p=1+\romega/m'>\romega$  and $q=1+\rgamma/\mu'>\rgamma$.

Moreover, if
\begin{equation}\label{1.14}
  p\le 1+\romega/2\,\text{when $N\ge 5$,\,$\gamma>0$,}\,\,\,\text{and}\,
  q\le 1+\rgamma/2\,\text{when $N\ge 6$, \, $\delta>0$,}
\end{equation}
then the previously found  solution is unique.

Beside the local theory described above, in \cite[Theorem~1.5 and Remarks~1.1, 1.3]{Dresda2} also existence of global solutions for arbitrary initial data was proved, provided the  parameters satisfy a further restriction.
\begin{theorem}[\bf Global existence]\label{theorem1.1} Let \eqref{1.3},\eqref{1.6},\eqref{1.11} hold and
\begin{equation}\label{1.15}
 p\le
\begin{cases}
2&\text{if $\gamma>0$, $\alpha=0$,}\\
\overline{m}&\text{if $\gamma>0$, $\alpha>0$,}
\end{cases}\qquad
q\le
\begin{cases}
2&\text{if $\delta>0$, $\beta=0$,}\\
\overline{\mu}&\text{if $\delta>0$, $\beta>0$.}
\end{cases}
\end{equation}
Then for any $U_0\in\cal{H}$ such that $u_0\in L^p(\Omega)$ when $\gamma>0$ and $p>\romega$,
${u_0}_{|\Gamma}\in L^q(\Gamma_1)$ when $\delta>0$ and $q>\rgamma$, problem \eqref{1.2} has a global weak solution, which is unique when also \eqref{1.14} holds.
\end{theorem}

\begin{figure}
\includegraphics[width=12.5truecm]{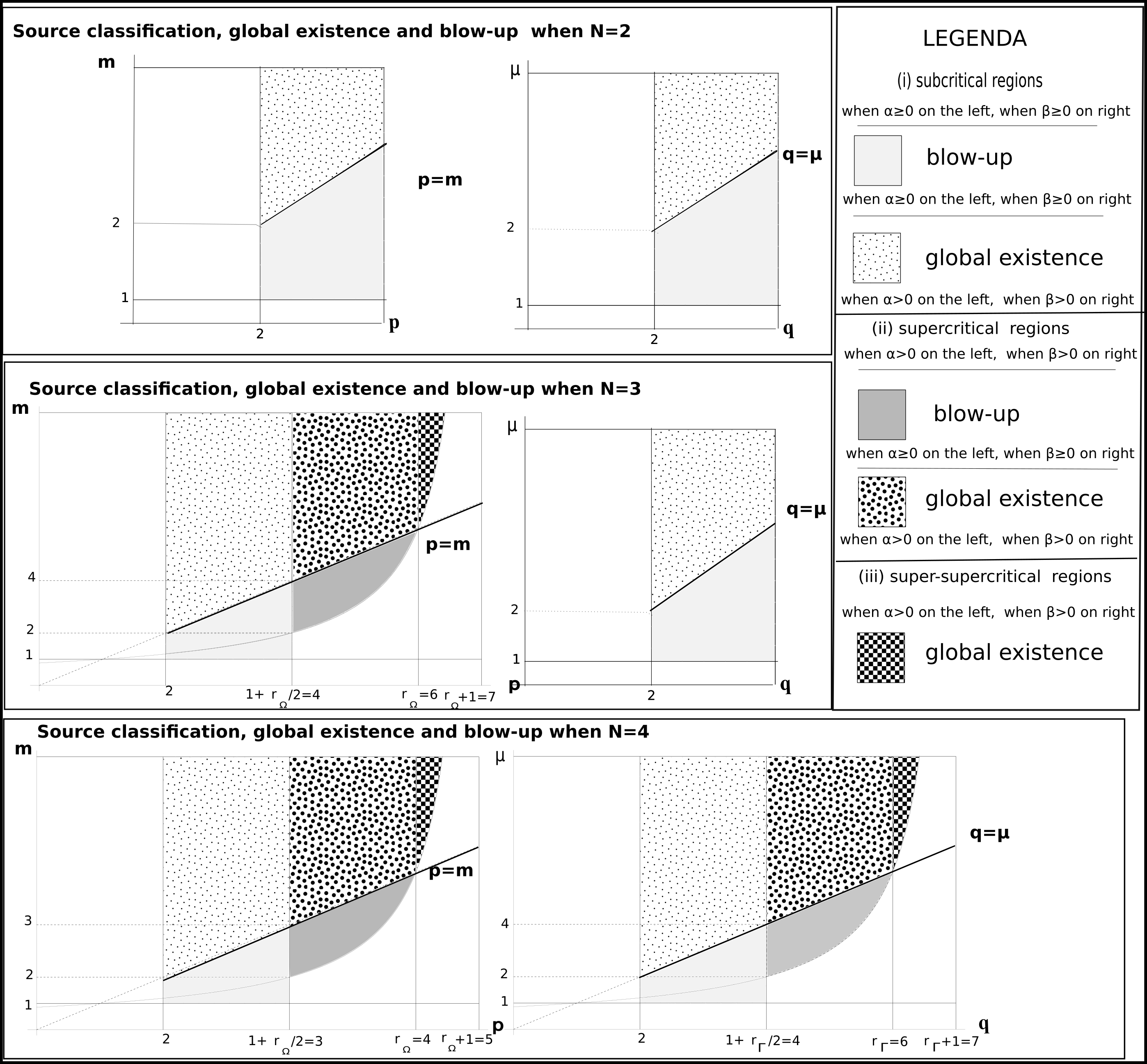}
\caption{The regions covered by \eqref{1.6}, and their subregions in which \eqref{1.15} holds or not, in dimensions $N=2,3,4$.}
\end{figure}

\begin{remark}\label{remark1.1}The sets in the planes $(p,m)$ and $(q,\mu)$, for which \eqref{1.6}  and \eqref{1.15} hold, and those for which \eqref{1.6} holds while \eqref{1.15} does not, depending on
the vanishing of $\alpha$ and $\beta$, are illustrated in
dimensions $N=2,3,4$ in Figure 1. As shown (when  $N\le 4$) in it, when \eqref{1.6} holds and one has
\begin{equation}\label{1.16}
\gamma>0\quad\text{and}\quad  p>
\begin{cases}
2&\text{if $\alpha=0$,}\\
\overline{m}&\text{if $\alpha>0$,}
\end{cases}
\end{equation}
i.e. when the first half of assumption \eqref{1.15} does not hold, one necessarily has $\overline{m},p<\romega$.
Indeed, if $\alpha=0$ then $p\le 1+\romega/2<\romega$ by \eqref{1.6}, and we can freely choose $m=2$,  while if $\alpha>0$ then by \eqref{1.6} and \eqref{1.16} one has $\overline{m}<1+\romega/\overline{m}'$, i.e. $\overline{m}^2-(\romega+1)\overline{m}+\romega>$, so $\overline{m}<\romega$ and a further application of \eqref{1.6} yields $p\le 1+\romega/\overline{m}'<1+\romega/\romega'=\romega$.
The same arguments show that then \eqref{1.6} holds and one has
\begin{equation}\label{1.17}
\delta>0\quad\text{and}\quad q>
\begin{cases}
2&\text{if $\beta=0$,}\\
\overline{\mu}&\text{if $\beta>0$,}
\end{cases}
\end{equation}
one necessarily has $\overline{\mu},q<\rgamma$.
\end{remark}
The optimality of assumption \eqref{1.15} was already discussed in \cite{Dresda1} when both damping terms  are linear, i.e. when
\begin{equation}\label{1.13}
a=b=0,\qquad  \alpha=0\quad\text{or}\quad m=2, \qquad \beta=0\quad\text{or}\quad \mu=2.
\end{equation}
 In this case, which includes the sourceless class, by \eqref{1.6} both sources are subcritical, so \eqref{1.14} holds, and assumption \eqref{1.15} trivializes to
\begin{equation}\label{1.18}
   p=2 \quad\text{when $\gamma>0$,} \qquad q=2 \quad\text{when $\delta>0$.}
\end{equation}
As a particular case of \cite[Theorem~1.5]{Dresda1} the following result holds.
\begin{theorem}[\bf Blow--up for linear damping]\label{theorem1.2}Let \eqref{1.3}, \eqref{1.6}, \eqref{1.13} hold, and
\begin{equation}\label{1.19}
  (\gamma,\delta)\not=(0,0),\qquad p>2 \quad\text{when $\gamma>0$,}\qquad q>2 \quad\text{when $\delta>0$.}
\end{equation}
Then, for any $U_0\in\cal{H}$ such that
\begin{multline}\label{1.25}
\cal{E}(U_0):=\tfrac 12 \|u_1\|_{H^0}^2+\tfrac 12 \int_\Omega |\nabla u_0|^2 dx+\tfrac 12 \int_{\Gamma_1} |\nabla_\Gamma u_0|_\Gamma^2d\sigma\\
-\tfrac\gamma p\int_\Omega |u_0|^p dx-\tfrac\delta q\!\int_{\Gamma_1} |u_0|^q d\sigma<0,
\end{multline}
 the unique maximal weak solution of \eqref{1.2} has $T_{\text{max}}<\infty$. Moreover
\begin{equation}\label{1.20}
  \lim_{t\to T_{\text{max}}^-}\|U(t)\|_{\cal{H}}=  \lim_{t\to T_{\text{max}}^-}\int_\Omega |u(t)|^p dx+\int_{\Gamma_1}|u(t)|^q d\sigma=\infty,
\end{equation}
where we can take $p=2$ when $\gamma=0$ and $q=2$ when $\delta=0$.
\end{theorem}

The relations between the parameters ranges \eqref{1.18} and \eqref{1.19}, which respectively yield global existence for (almost) all data and blow--up for suitable data (which trivially exist), is clearest when separately considering the previously introduced model classes A--D). This comparison is made explicit, for the readers' convenience, in Table~1. In it
$\checkmark$ stands for no assumptions and $\cal{X}$ for class exclusion. 

\bigskip
\begin{center}
\renewcommand{\arraystretch}{1.3}
\setlength{\arrayrulewidth}{0.6pt}
{\sc Table 1: \rm\eqref{1.18} vs. \eqref{1.19} for the model classes A--D).}\\ \smallskip
\begin{tabular}{|c|c|c|c|c|}
\hline
Case \eqref{1.17}  &A) $\gamma=0$, $\delta=0$  &B) $\gamma=0$, $\delta>0$ &C) $\gamma>0$, $\delta=0$ & D) $\gamma>0$, $\delta>0$ \\
\hline
\eqref{1.18}&$\checkmark$&$q=2$ &$p=2$ & $p=2$, $q=2$ \\
\hline
\eqref{1.19}&$\cal{X}$ &$q>2$ &$p>2$ & $p>2$, $q>2$ \\
\hline
\end{tabular}
\end{center}
\bigskip
Theorem~\ref{theorem1.1}, and consequently also Theorem~\ref{theorem1.2}, is sharp in the classes A--C), while in the class D) the combined answer given by them is incomplete. Indeed, when $p=2$, $q>2$ and  when $p>2$, $q=2$ no information is given. This easy case exhibits two difficulties in the analysis of the class D) which will persists in the general case:

\begin{itemize}
  \item[$\bullet$)] even if a source is superlinear, the linearity of the other one may inhibit global nonexistence arguments,
  \item[$\bullet\bullet$)]  when the growth of one source dominated the growth of the corresponding damping term, but the opposite domination holds for the other couple, the solutions behavior may remain undetermined.
\end{itemize}

\subsection{Main results} We shall present our main results for \eqref{1.2}  by distinguishing among the previously introduces model classes A--D). Clearly when \eqref{1.2} is sourceless Theorem~\ref{theorem1.2} assures global existence for (almost) all data, so the class A) needs no further attention.

We start by making some remarks on class B), when only a boundary source is present.  This class is covered by Theorem~\ref{theorem1.2} when the damping is linear, but the situation is quite different in the nonlinear case. Indeed it is possible to find some blow--up results in the literature concerning class B) for the related family \eqref{1.4}, such as  \cite{Ha1, Ha2, ZhangHu}. Unfortunately all the proofs of this type of results are, to the author's knowledge,  adaptations of the classical arguments in \cite{georgiev,levserr} and are, at some point, problematic. Referring to the quoted papers, in \cite[p. 868]{ZhangHu} the authors treats the norms $\|\cdot\|_{L^p(\Omega)}$ and $\|\cdot\|_{L^p(\Gamma_1)}$ as equivalent, while in \cite[p. 333]{Ha1} (the same argument being used in \cite{Ha2}) the author implicitly uses boundedness of the function $t\mapsto\|u_t(t)\|_{L^2(\Omega)}$ while proving the finite time blow--up of an auxiliary functional, which in turns yields finite time blow--up of $t\mapsto\|u_t(t)\|_{L^2(\Omega)}$.
To the author's understanding the arguments in \cite{georgiev,levserr} cannot be adapted to wave equation with boundary nonlinear damping and sources and, since this is the state of the arts, is still a challenging  open problem to prove blow--up results in class B) for problem \eqref{1.4}. In the present paper we shall not  deal with this class for \eqref{1.2}.

Our first main result concerns the model class C).
\begin{theorem}[\bf Global nonexistence and blow--up with interior source]\label{theorem1.3}\phantom{.}\\
Let \eqref{1.3}, \eqref{1.6} hold, $\gamma>0$, $\delta=0$ and
\begin{equation}\label{1.21}
p>\begin{cases}
2&\text{if $\alpha=0$,}\\
\overline{m}&\text{if $\alpha>0$,}
\end{cases}
\qquad
\overline{\mu}<1+p/2\quad\text{when $\beta>0$.}
\end{equation}
Then, for any $U_0\in\cal{H}$ such that
\begin{equation}\label{1.22}
\cal{E}(U_0)=\tfrac 12 \|u_1\|_{H^0}^2+\tfrac 12 \int_\Omega |\nabla u_0|^2 dx+\tfrac 12 \int_{\Gamma_1} |\nabla_\Gamma u_0|_\Gamma^2 d\sigma-\tfrac\gamma p\int_\Omega |u_0|^p dx<0,
\end{equation}
and any maximal weak solution  of \eqref{1.2}   in $[0,T_{\text{max}})$,  one has
$T_{\text{max}}<\infty$ and
\begin{equation}\label{1.23}
  \varlimsup_{t\to T_{\text{max}}^-}\|U(t)\|_{\cal{H}}=  \varlimsup_{t\to T_{\text{max}}^-}\int_\Omega |u(t)|^p dx+\int_{\Gamma_1}|u(t)|^2 d\sigma=\infty.
\end{equation}
Finally, when $N\le 4$ or $p\le 1+\romega/2$, we can replace $\varlimsup\limits_{t\to T_{\text{max}}^-}$ with $\lim\limits_{t\to T_{\text{max}}^-}$ in \eqref{1.23}.
\end{theorem}

\begin{remark}\label{remark1.3} As already pointed out in Remark~\ref{remark1.1}, by \eqref{1.21} we necessarily have $\overline{m},p<\romega$. Since $\delta=0$, in Theorem~\ref{theorem1.3} super--supercritical sources are not considered, so ruling out conditions \eqref{1.11}--\eqref{1.12} from local existence theory.
We also point out that \eqref{1.23} holds for \emph{all possible} maximal weak solutions of \eqref{1.2}, also when uniqueness is unknown.
\end{remark}
The relation between the parameter ranges \eqref{1.15} and \eqref{1.21}, which respectively yield global existence for (almost) all data and blow--up for suitable data (which trivially exist) is clearest when separately considering the model problems Ca--Cd),  as we do in Table~2 for the reader's convenience. 

\bigskip
\begin{center}
\renewcommand{\arraystretch}{1.3}
\setlength{\arrayrulewidth}{0.6pt}
{\sc Table 2: \rm\eqref{1.15} vs. \eqref{1.21} for the model problems Ca--Cd).}\\ \smallskip
\begin{tabular}{|c|c|c|c|c|}
\hline
 $\gamma>0=\delta$ &$\alpha=\beta=0$  & $\alpha=0<\beta$ & $\alpha>0=\beta$ &  $\alpha,\beta>0$ \\
\hline
\eqref{1.15}&$p=2$&$p=2$ &$p\le\overline{m}$ & $p\le\overline{m}$ \\
\hline
\eqref{1.21}&$p>2$ &$p>2$, $\overline{\mu}<1+p/2$   &$p>\overline{m}$ & $p>\overline{m}$, $\overline{\mu}<1+p/2$\\
\hline
\end{tabular}
\end{center}

\bigskip

It makes clear that Theorem~\ref{theorem1.1}, and consequently also Theorem~\ref{theorem1.3}, is sharp for the model problems Ca) and Cc), i.e. when boundary source and damping are not present in \eqref{1.2}. By the contrary, when dealing with the model problems Cb) and Cd) and
\begin{equation}\label{1.24}
p>\begin{cases}
2&\text{if $\alpha=0$,}\\
\overline{m}&\text{if $\alpha>0$,}
\end{cases}
\qquad
\overline{\mu}\ge 1+p/2,
\end{equation}
Theorems~\ref{theorem1.1} and \ref{theorem1.3} give no information. When \eqref{1.24} holds the growth of the interior source dominates the one of the corresponding damping term, while the boundary damping term  has no a homologous counterpart and, consequently, can be controlled only by the \emph{transversal} influence of the interior source. This type of transversal control was already pointed out for the class C) of the related family \eqref{1.4} (without the term $u$ on $\Gamma_1$) in \cite{stable} and  subsequently improved in \cite{fisvit2}, where the exact assumption $\overline{\mu}<1+p/2$ appears. As to the author's knowledge such an assumption has been skipped only when $N=1$ and $\Omega$ is a suitably large interval (see \cite{FLZ}).

Our second main result concerns the model class D).
\begin{theorem}[\bf Global nonexistence and blow--up with two sources]\label{theorem1.4}\phantom{.}\\
Let \eqref{1.3}, \eqref{1.6} hold, $\gamma, \delta>0$ and
\begin{equation}\label{1.24bis}
p>\begin{cases}
2&\text{if $\alpha=0$,}\\
\overline{m}&\text{if $\alpha>0$,}
\end{cases}
\qquad
q>2,\quad \overline{\mu}<\max\{q,1+p/2\}\quad\text{when $\beta>0$.}
\end{equation}
Then, for any $U_0\in\cal{H}$ such that \eqref{1.25} holds
and any maximal weak solution \begin{footnote}{at least one of them exists}\end{footnote} of \eqref{1.2} in $[0,T_{\text{max}})$, one has
$T_{\text{max}}<\infty$ and
\begin{equation}\label{1.26}
  \varlimsup_{t\to T_{\text{max}}^-}\|U(t)\|_{\cal{H}}=  \varlimsup_{t\to T_{\text{max}}^-}\int_\Omega |u(t)|^p dx+\int_{\Gamma_1}|u(t)|^q d\sigma=\infty.
\end{equation}
Finally, when $N\le 4$, or $N=5$ and $p\le 1+\romega/2=8/3$, or $N\ge 6$, $p\le 1+\romega/2$, $q\le 1+\rgamma/2$, we can replace $\varlimsup_{t\to T_{\text{max}}^-}$ with $\lim_{t\to T_{\text{max}}^-}$ in \eqref{1.26}.
\end{theorem}
\begin{remark}\label{remark1.5} As already pointed out in Remark~\ref{remark1.1}, by \eqref{1.6} and \eqref{1.24bis} we necessarily have $\overline{m},p<\romega$.
By \eqref{1.6} and \eqref{1.24bis} it also follows that $\overline{\mu},q<\rgamma$. Indeed, when $\beta=0$, by \eqref{1.6}
we have $q\le 1+\rgamma/2<\rgamma$ and we can freely choose $\mu=2$, while when $\beta>0$ by \eqref{1.24bis} either $\overline{\mu}<q$ or $\overline{\mu}<1+p/2$. In the first case, see also Figure~1, by \eqref{1.6} one has  $\overline{\mu}<1+\rgamma/\overline{\mu}'$, i.e. $\overline{\mu}^2-(\rgamma+1)\overline{\mu}+\rgamma<0$, so $\overline{\mu}<\rgamma$ and consequently, using \eqref{1.6} again, $q\le 1+\rgamma/\overline{\mu}'<1+\rgamma/\rgamma'=\rgamma$. In the second case, since $p<\romega$, we have $\overline{\mu}<1+\romega/2$. But, by \eqref{1.5}, $1+\romega/2\le \rgamma$ for all $N\ge 2$, so $\overline{\mu}<\rgamma$ and, as in the previous case, $q<\rgamma$.
Hence, also in Theorem~\ref{theorem1.4} super--supercritical sources are not considered, so ruling out conditions \eqref{1.11}--\eqref{1.12} from local existence theory. Also in this case \eqref{1.26} holds for \emph{all possible} maximal weak solutions of \eqref{1.2}, also when uniqueness is unknown.
\end{remark}
Clearly \eqref{1.24bis} does noth exhaust all possible parameters values for which \eqref{1.15} does not hold in class D). Also for this class it is useful to separately considering the model problems Da--Dd). This comparison is made, for the reader's convenience, in Table~3. 

\bigskip

\begin{center}
\renewcommand{\arraystretch}{1.3}
\setlength{\arrayrulewidth}{0.6pt}
{\sc Table 3: \rm\eqref{1.15} vs \eqref{1.24bis} for the model problems Da--Dd).}\\ \smallskip
\begin{tabular}{|c|c|c|c|c|}
\hline
 $\gamma,\delta>0$ &$\alpha=\beta=0$  & $\alpha=0<\beta$ & $\alpha>0=\beta$ &  $\alpha,\beta>0$ \\
\hline
\eqref{1.15}&$p=2$,  &$p=2$,                &$p\le\overline{m}$,   & $p\le\overline{m}$,  \\
            &$q=2$   &$q\le \overline{\mu}$ &$q=2$                 & $q\le\overline{\mu}$ \\
\hline
\eqref{1.24bis}&$p>2$,  &$p>2$,   &$p>\overline{m}$, & $p>\overline{m}$, \\
    &$q>2$              &$q>\begin{cases}2\,&\text{if $\overline{\mu}<1+\frac p2$}\\ \overline{\mu}\,&\text{if $\overline{\mu}\ge 1+\frac p2$}\end{cases}$   &$q>2$ & $q>\begin{cases}2\,&\text{if $\overline{\mu}<1+\frac p2$}\\ \overline{\mu}\,&\text{if $\overline{\mu}\ge 1+\frac p2$}\end{cases}$\\

\hline
\end{tabular}
\end{center}

\bigskip
It clearly shows the following facts.
In models Da) and Dc) Theorem~\ref{theorem1.4} predicts finite time blow--up of solutions only when \emph{both} inequalities in \eqref{1.15} do not hold. So no information is given in the following two cases:
\renewcommand{\labelenumi}{\roman{enumi})}
\begin{enumerate}
\item when the linear boundary source inhibits global nonexistence arguments, i.e. when
$$p>\begin{cases}
2&\text{if $\alpha=0$,}\\
\overline{m}&\text{if $\alpha>0$,}
\end{cases}
\qquad
q=2;$$
\item when the growth of the interior damping dominates the one of the interior source, but the boundary source has no a damping homologous counterpart, i.e. when
$$p\le\begin{cases}
2&\text{if $\alpha=0$,}\\
\overline{m}&\text{if $\alpha>0$,}
\end{cases}
\qquad
q>2.$$
\end{enumerate}
These two cases exactly correspond to the  difficulties $\bullet)$ and $\bullet\bullet)$ ( which show up at their top level for these models) already emphasized for class D) when damping terms are linear.

In model problems Db) and Dd)  Theorem~\ref{theorem1.4} predicts blow--up of solutions when both inequalities in \eqref{1.15} do not hold, that is  when \eqref{1.16}--\eqref{1.17} hold, but \emph{not only} in this case. Indeed it also partially cover the case in which \eqref{1.16} holds while \eqref{1.17} does not, this part exactly corresponding to the one covered in models Cb) and Cd).

It is worth observing that the parameter restriction \eqref{1.24bis} in the literature concerning class D) for the family \eqref{1.4}, see \cite{BLWarsaw,BRT, BLHindawi}, represents the state of the art.

Finally we remark that \eqref{1.24bis} is the exact complementary of \eqref{1.15}, so Theorems~\ref{theorem1.3}--\ref{theorem1.4} are both sharp, in  the "diagonal" case $\alpha=\beta$, $\gamma=\delta$, $m=\mu$ and $p=q$.

Theorems~\ref{theorem1.3}--\ref{theorem1.4} are suitably recombined particular cases of our more general blow--up results Theorems~\ref{theorem5.3} and \ref{theorem5.4}, which are essentially based on our main global nonexistence results Theorem~\ref{theorem3.1} and \ref{theorem4.1}. Their proofs both rely on suitable non--trivial adaptations of the techniques in \cite{BLWarsaw,BRT, fisvit2,BLHindawi,georgiev,levserr, stable}.

\subsection{Organization of the paper}
The sequel of the paper is organized as follows:
\begin{enumerate}
  \item in Section~\ref{section 2}  we give some background on the functional spaces used and  on a linear version of problem \eqref{1.1};
  \item Section~\ref{section 3} is devoted to give our main assumptions, to introduce weak solutions of \eqref{1.1}  and to some preliminary results;
  \item in Section~\ref{section4}  we state and prove  our first main global nonexistence result for problem \eqref{1.1}, dealing with the case when two sources are present in it;
  \item in Section~\ref{section5}  we state and prove our second main global nonexistence result for problem \eqref{1.1}, dealing with the case in which $g$ may vanish;
  \item Section~\ref{section6} is devoted to recall the local theory from \cite{Dresda2}, to  give our main blow--up results and to show how Theorems~\ref{theorem1.3}--\ref{theorem1.4} follow from them.
\end{enumerate}

\section{Background}\label{section 2}
\subsection{Notation}
We shall adopt the standard notation for
(real) Lebesgue and Sobolev spaces in $\Omega$ (see \cite{adams})
and $\Gamma$ (see \cite{grisvard}). Moreover $\|\cdot\|_\rho:=\|\cdot\|_{L^\rho(\Omega)}$
and $\|\cdot\|_{\rho,\Gamma'}:=\|\cdot\|_{L^\rho(\Gamma')}$ for $\rho\in [1,\infty)$ and $\Gamma'\subseteq\Gamma$ measurable.

Given any Banach space $X$ we shall denote by $\langle\cdot,\cdot\rangle_X$ the duality product between $X$ and its dual $X'$, and we shall use the standard notation for $X$--valued Bochner--Lebesgue and
Bochner--Sobolev spaces in a real interval. Moreover $\cal{L}(X,Y)$ will denote the class of linear bounded operators from $X$ to another Banach space $Y$.

 Given $\alpha\in
L^\infty(\Omega)$, $\beta\in L^\infty(\Gamma_1)$, $\alpha,\beta\ge 0$, $-\infty\le c<d\le\infty$
and $\rho\in [1,\infty)$ we shall respectively denote by
$\lambda_\alpha$, $\lambda'_\alpha$, $\lambda_\beta$, $\lambda'_\beta$ the measures respectively defined
in $\Omega$, $\R\times\Omega$, $\Gamma_1$, $\R\times\Gamma_1$,  by $d\lambda_\alpha=\alpha\,dx$,
$d\lambda'_\alpha=\alpha\,dt\,dx$, $d\lambda_\beta=\beta\,d\sigma$, $d\lambda'_\beta=\beta\,dt\,d\sigma$, and by $L^\rho(\Omega;\lambda_\alpha)$,  $L^\rho((c,d)\times\Omega;\lambda'_\alpha)$, $L^\rho(\Gamma_1;\lambda_\beta)$,  $L^\rho((c,d)\times\Gamma_1;\lambda'_\beta)$ the corresponding Lebesgue spaces.
The equivalence classes with respect to $\lambda_\alpha$ and $\lambda'_\alpha$ a.e. equivalences will be denoted by $[\cdot]_\alpha$, those with respect to $\lambda_\beta$ and $\lambda'_\beta$ equivalences  by $[\cdot]_\beta$.
By the density of $C_c((c,d)\times\Omega)$ and  $C_c((c,d)\times\Gamma_1)$, respectively in $L^\rho((c,d)\times\Omega;\lambda'_\alpha)$ and $L^\rho((c,d)\times\Gamma_1;\lambda'_\beta)$, see \cite[Theorem~2.18 p. 48 and Theorem~3.14 p. 68]{rudin}, one can prove, as in the standard case, that
\begin{equation}\label{2.1}
\begin{gathered}
 L^\rho(c,d;L^\rho(\Omega,\lambda_\alpha))\simeq L^\rho((c,d)\times\Omega;\lambda'_\alpha),\\
 L^\rho(c,d;L^\rho(\Gamma_1,\lambda_\beta))\simeq L^\rho((c,d)\times\Gamma_1;\lambda'_\beta).
 \end{gathered}
\end{equation}
We recall some well--known preliminaries on the Riemannian gradient
on the $C^1$ compact manifold $\Gamma$, referring to
\cite{taylor} for more details and proofs in the smooth
setting and to \cite{sternberg} for a general background on
differential geometry on $C^1$  manifolds. The interest reader may also see \cite{mugnvit} in the non--compact case. We  denote
by $(\cdot,\cdot)_\Gamma$ the metric inherited from $\R^N$,
  given in local coordinates
$(y_1,\ldots,y_{N-1})$ by $(g_{ij})_{i,j=1,\ldots,N-1}$, and
$|\cdot|_\Gamma^2=(\cdot,\cdot)_\Gamma$. We denote by $d\sigma$ the natural
volume element on $\Gamma$, given by $\sqrt{\operatorname{det} (g_{ij})}
\,\,dy_1\wedge\ldots\wedge dy_{N-1}$. The Riemannian
gradient is given in local coordinates by $\nabla_\Gamma
u=g^{ij}\,\partial_j u\,\,\partial_i$ for any $u\in H^1(\Gamma)$, where $(g^{ij})=(g_{ij})^{-1}$.
It is well--known, see
\cite{quarteroni,mugnvit,taylor} that the norm
$\|u\|_{H^1(\Gamma)}^2=\|u\|_{2,\Gamma}^2+\|\nabla_\Gamma
u\|_{2,\Gamma}^2$,  where $\|\nabla_\Gamma
u\|_{2,\Gamma}^2:=\int_\Gamma |\nabla_\Gamma u|_\Gamma^2$, is
equivalent in $H^1(\Gamma)$  to the standard one.

\subsection{Functional setting}
Given $\alpha\in L^\infty(\Omega)$, $\beta\in L^\infty(\Gamma_1)$, $\alpha,\beta\ge 0$
and $\rho\in [2,\infty)$ we recall the reflexive spaces introduced in \cite{Dresda1}, that is
\begin{align*}
&L^{2,\rho}_\alpha(\Omega)=\{u\in L^2(\Omega): \alpha^{1/\rho}u\in
L^\rho(\Omega)\}, &&\|\cdot\|_{L^{2,\rho}_\alpha(\Omega)}=\|\cdot\|_2+\|\alpha^{1/\rho}
\cdot\|_\rho,\\
& L^{2,\rho}_\beta(\Gamma_1)=\{u\in L^2(\Gamma_1):
\beta^{1/\rho}u\in L^\rho(\Gamma_1)\},
&&\|\cdot\|_{L^{2,\rho}_\beta(\Gamma_1)}=\|\cdot\|_{2,\Gamma_1}+\|\beta^{1/\rho}
\cdot\|_{\rho,\Gamma_1},
\end{align*}
as well as the trivial embeddings  and boundedness properties
\begin{equation}\label{2.3}
\begin{aligned}
L^{2,\rho}_\alpha(\Omega)\hookrightarrow L^2(\Omega),\qquad &[\cdot]_\alpha\in
\cal{L}(L^{2,\rho}_\alpha(\Omega),L^\rho(\Omega,\lambda_\alpha))\\
L^{2,\rho}_\beta(\Gamma_1)\hookrightarrow L^2(\Gamma_1),\qquad &[\cdot]_\beta\in
\cal{L}(L^{2,\rho}_\beta(\Gamma_1),L^\rho(\Gamma_1,\lambda_\beta)).
\end{aligned}
\end{equation}
All operators in \eqref{2.3} trivially have dense range (see \cite{Dresda1}), so by applying \cite[Theorem~5.11--3, Chapter 5, p. 280]{Ciarlet}, or \cite[Corollary~2.18 p. 45]{brezis2}, we have the embeddings
\begin{equation}\label{2.4}
\begin{aligned}
[L^2(\Omega)]'\hookrightarrow [L^{2,\rho}_\alpha(\Omega)]',\qquad &[L^\rho(\Omega,\lambda_\alpha)]'\hookrightarrow [L^{2,\rho}_\alpha(\Omega)]'\\
[L^2(\Gamma_1)]'\hookrightarrow [L^{2,\rho}_\beta(\Omega)]',\qquad &[L^\rho(\Gamma_1,\lambda_\beta)]'\hookrightarrow [L^{2,\rho}_\beta(\Gamma_1)]'.
\end{aligned}
\end{equation}
As usual we shall identify $[L^2(\Omega)]'\simeq  L^2(\Omega)$ and $[L^2(\Gamma_1)]'\simeq  L^2(\Gamma_1)$. These identifications, essentially made in the distribution sense, make impossible to identify
$[L^\rho(\Omega,\lambda_\alpha)]'$ with $L^{\rho'}(\Omega,\lambda_\alpha)$, the same remark applying to measures $\lambda'_\alpha$, $\lambda_\beta$ and  $\lambda'_\beta$.
We shall identify all spaces in \eqref{2.4} with the corresponding subspaces
of $[L^{2,\rho}_\alpha(\Omega)]'$ or $[L^{2,\rho}_\beta(\Gamma_1)]'$.

For any $\xi \in L^{\rho'}(\Omega,\lambda_\alpha)$, even if $\xi$ is not a.e. well-defined in $\Omega$, since it takes arbitrary values in the possibly large set where $\alpha$ vanishes, the function $\alpha\xi$ is well--defined a.e. in it, and actually $\alpha\xi\in L^1(\Omega)$. Moreover, by the form of the Riesz isomorphism in $L^\rho(\Omega,\lambda_\alpha)$,  we can represent   $[L^\rho(\Omega,\lambda_\alpha)]'$ as $\{\alpha\xi, \xi \in L^{\rho'}(\Omega, \lambda_\alpha)\}$.
Using the same arguments on $\Gamma_1$ we have that for any $\eta\in L^{\rho'}(\Gamma_1,\lambda_\beta)$ we have $\beta\eta\in L^1(\Gamma_1)$ and the following  identifications hold
\begin{equation}\label{identificazioni0}
[L^\rho(\Omega,\lambda_\alpha)]'\simeq\{\alpha\xi, \xi \in L^{\rho'}(\Omega, \lambda_\alpha)\},
\quad
[L^\rho(\Gamma_1,\lambda_\beta)]'\simeq\{\beta\eta, \eta \in L^{\rho'}(\Gamma_1, \lambda_\beta)\}.
\end{equation}
In the sequel we shall also use , for any $\alpha\in
L^\infty(\Omega)$, $\beta\in L^\infty(\Gamma_1)$, $\alpha,\beta\ge 0$, $-\infty\le c<d\le\infty$
and $\rho\in [2,\infty)$, the spaces $L^\rho(c,d;L^{2,\rho}_\alpha (\Omega))$ and  $L^\rho(c,d;L^{2,\rho}_\beta (\Gamma_1))$. Trivially, by \eqref{2.1}, \eqref{2.3} and \eqref{2.4},
\begin{equation}\label{2.5}
\begin{aligned}
&[\cdot]_\alpha\in \cal{L}(L^\rho(c,d;L^{2,\rho}_\alpha(\Omega)),L^\rho(c,d;L^\rho(\Omega,\lambda_\alpha)))\\
&[\cdot]_\beta\in \cal{L}(L^\rho(c,d;L^{2,\rho}_\beta(\Gamma_1)),L^\rho(c,d;L^\rho(\Gamma_1,\lambda_\beta)))\\
&L^{\rho'}((c,d)\times\Omega,\lambda'_\alpha)\hookrightarrow L^{\rho'}(c,d;[L^{2,\rho}_\alpha(\Omega)]')\\
&L^{\rho'}((c,d)\times\Gamma_1,\lambda'_\beta)\hookrightarrow L^{\rho'}(c,d;[L^{2,\rho}_\beta(\Gamma_1)]').
\end{aligned}
\end{equation}
In the sequel we shall treat the embeddings in last two lines of \eqref{2.5} as identifications.
The same arguments used before to get \eqref{identificazioni0} allows us to make the further identifications
\begin{equation}\label{identificazioni}
\begin{aligned}
[L^\rho(c,d;L^\rho(\Omega,\lambda_\alpha)]'\simeq &\{\alpha\xi,\quad \xi \in L^{\rho'}(c,d;L^{\rho'}(\Omega, \lambda_\alpha))\},\\
[L^\rho(c,d;L^\rho(\Gamma_1,\lambda_\beta)]'\simeq & \{\beta\eta, \quad \eta \in L^{\rho'}(c,d;L^{\rho'}(\Gamma_1, \lambda_\beta))\},
\end{aligned}
\end{equation}
where, when $-\infty<c<d<\infty$,  $\alpha\xi\in L^1((c,d)\times\Omega)$ and $\beta\eta\in L^1((c,d)\times\Gamma_1)$.

Next, given $\rho,\theta\in [2,\infty)$ and
$-\infty\le c<d\le\infty$ we introduce the  space
\begin{equation}\label{2.6}
L^{2,\rho,\theta}_{\alpha,\beta}(c,d)=L^\rho(c,d\, ;
L^{2,\rho}_\alpha (\Omega))\times L^\theta(c,d\, ;
L^{2,\theta}_\beta (\Gamma_1)),
\end{equation}
together with its right--local version
\begin{equation}\label{2.6bis}
L^{2,\rho,\theta}_{\alpha,\beta,\loc}([c,d))=L^\rho_\loc([c,d))\, ;
L^{2,\rho}_\alpha (\Omega))\times L^\theta_\loc([c,d))\, ;
L^{2,\theta}_\beta (\Gamma_1)).
\end{equation}

We respectively endow the Hilbert spaces $H^0$  and $H^1$ introduced in \S~\ref{subsection1.3} with the standard inner product
given, for $w_i=(u_i,v_i)\in H^0$, $i=1,2$, by
\begin{equation}\label{2.7}
  (w_1,w_2)_{H^0}=\int_\Omega u_1u_2\,dx +\int_{\Gamma_1} v_1v_2\,d\sigma,
\end{equation}
and with the inner product
\begin{equation}\label{2.8}
  (u,v)_{H^1}=\int_\Omega \nabla u\nabla v \,dx+\int_{\Gamma_1} (\nabla_\Gamma u,\nabla_\Gamma v)_\Gamma \,d\sigma+\int_{\Gamma_1} u v \,d\sigma,\quad u,v\in H^1.
\end{equation}
Its associated norm $\|\cdot\|_{H^1}=(\cdot,\cdot)_{H^1}^{1/2}$ is equivalent to the standard one inherited from the product.
 We also introduce, for any $\alpha\in
L^\infty(\Omega)$, $\beta\in L^\infty(\Gamma_1)$, $\alpha,\beta\ge 0$,
$\rho,\theta\in [2,\infty)$, the Banach
space
\begin{equation}\label{2.9}
H^{1,\rho,\theta}_{\alpha,\beta}=H^1\cap
[L^{2,\rho}_\alpha(\Omega)\times L^{2,\theta}_\beta(\Gamma_1)],\quad \|\cdot\|_{H^{1,\rho,\theta}_{\alpha,\beta}}=\|\cdot\|_{H^1}+\|\cdot\|_{L^{2,\rho}_\alpha(\Omega) \times L^{2,\rho}_\beta(\Gamma_1)}
\end{equation}
and
\begin{equation}\label{2.10}
 H^{1,\rho,\theta}=H^{1,\rho,\theta}_{1,1}=H^1\cap
[L^\rho(\Omega)\times L^\theta(\Gamma_1)].
\end{equation}
Trivially
\begin{equation}\label{2.11}
H^{1,\rho,\theta}\hookrightarrow H^{1,\rho,\theta}_{\alpha,\beta}\hookrightarrow H^1.
\end{equation}
Finally, beside the main phase space $\cal{H}$ introduced in \eqref{1.10}, we also introduce the auxiliary phase spaces
\begin{equation}\label{2.13}
 \cal{H}^{\rho,\theta}=H^{1,\rho,\theta}\times H^0,
\end{equation}
which trivially does not coincides with $\cal{H}$ only when $\rho>\romega$ or $\theta>\rgamma$.
Although in our main result we shall not consider super--supercritical sources, for which these  spaces are needed, we shall use them when introducing weak solutions of \eqref{1.1}. They can be useful in further studies.
\subsection{Weak solutions for a linear version of \eqref{1.1}}
 We consider the linear
evolution boundary value problem
\begin{equation}\label{2.14}
\begin{cases} u_{tt}-\Delta u=\xi \qquad &\text{in
$(0,T)\times\Omega$,}\\
u=0 &\text{on $(0,T)\times \Gamma_0$,}\\
u_{tt}+\partial_\nu u-\Delta_\Gamma u=\eta\qquad &\text{on
$(0,T)\times \Gamma_1$,}
\end{cases}
\end{equation}
where $0<T<\infty$ and $\xi=\xi(t,x)$, $\eta=\eta(t,x)$ are given
forcing terms of the form
\begin{equation}\label{2.15}
\left\{
\begin{alignedat}3
&\xi=\xi_1+\alpha \xi_2,\qquad && \xi_1\in L^1(0,T;L^2(\Omega)), \quad && \xi_2\in L^{\rho'}(0,T;L^{\rho'}(\Omega, \lambda_\alpha)), \\
&\eta=\eta_1+\beta \eta_2,\qquad && \eta_1\in
L^1(0,T;L^2(\Gamma_1)), \quad && \eta_2\in
L^{\theta'}(0,T;L^{\theta'}(\Gamma_1, \lambda_\beta)),
\end{alignedat}\right.
\end{equation}
where $\alpha\in L^\infty(\Omega)$, $\beta\in L^\infty(\Gamma_1)$,
$\alpha, \beta\ge 0$ and $\rho,\theta\in [2,\infty)$. Hence $\xi\in L^1((0,T)\times\Omega)$, $\eta\in\L^1((0,T)\times\Gamma_1)$ and,  by \eqref{identificazioni} and \eqref{2.6},
\begin{equation}\label{2.16}
\xi\in L^1(0,T; [L^{2,\rho}_\alpha(\Omega)]'),\qquad
\eta\in L^1(0,T; [L^{2,\theta}_\beta(\Gamma_1)]'),
\end{equation}
so that the following definition makes sense.
\begin{definition}\label{definition2.1}
Let $\xi$ and $\eta$ be given by \eqref{2.15}. By a {\em weak solution} of \eqref{2.14} in $[0,T]$   we mean
$u\in L^\infty(0,T;H^1)\cap W^{1,\infty}(0,T;H^0)$, $u'\in L^{2,\rho,\theta}_{\alpha,\beta}(0,T)$,
such that the distribution identity
\begin{multline}\label{2.17}
\int_0^T\left[-(u',\phi')_{H^0}+\int_\Omega \nabla u\nabla
\phi \,dx+\int_{\Gamma_1} (\nabla_\Gamma
u,\nabla_\Gamma\phi)_\Gamma d\sigma\right.\\\left.-\int_\Omega
\xi\phi \,dx-\int_{\Gamma_1}\eta\phi \,d\sigma\right]=0
\end{multline}
holds for all $\phi\in C_c((0,T);H^1)\cap C^1_c((0,T);H^0)\cap
L^{2,\rho,\theta}_{\alpha,\beta}(0,T)$.
\end{definition}
Clearly, for any weak solution $u$ of \eqref{2.14}, one has $u'=(u_t,(u_{|\Gamma})_t)$.
Since the two components of $u'$, respectively acting in $(0,T)\times\Omega$ and in $(0,T)\times\Gamma_1$,
cannot be confused, with a slight abuse we shall denote, for simplicity,  $(u_{|\Gamma})_t=u_t$.
Hence  we shall {\em
systematically denote in the paper}
\begin{equation}\label{2.18}
u'=(u_t,u_t)\qquad\text{and}\quad U=(u,u')\in
L^\infty(0,T;\cal{H}).
\end{equation}

We recall  \cite[Lemma 2.2]{Dresda1}.
\begin{lemma}\label{lemma2.1}
Any weak solution $u$ of
\eqref{2.14} enjoys the further regularity $U\in C([0,T];\cal{H})$. Moreover it
satisfies the following identities:
\begin{enumerate}
  \item the energy identity
$$\frac 12\|u'\|_{H^0}^2\!+\!\frac 12\!
\int_\Omega\!\! |\nabla u|^2 dx+\!\frac 12 \!\int_{\Gamma_1}\!\!|\nabla u|_\Gamma^2 d\sigma\bigg|_s^t\!\!=\!\!\int_s^t\!\!\left(\int_\Omega \xi
u_t \,dx+\!\!\int_{\Gamma_1}\eta u_t \,d\sigma\right)d\tau  $$ for $0\le s\le
t\le T$;
  \item the generalized distribution identity
\begin{align*}
(u',\phi)_{H^0}\Big|_0^T+
\int_0^T&\left[-(u',\phi')_{H^0}+\int_\Omega \nabla u\nabla
\phi \,dx\right. \\ &\left.+\int_{\Gamma_1} (\nabla_\Gamma
u,\nabla_\Gamma\phi)_\Gamma d\sigma-\int_\Omega
\xi\phi \,dx -\int_{\Gamma_1}\eta\phi \,d\sigma\right]=0
\end{align*}
 for all $\phi\in
C([0,T];H^1)\cap C^1([0,T];H^0)\cap
L^{2,\rho,\theta}_{\alpha,\beta}(0,T)$.
\end{enumerate}
\end{lemma}
\section{Preliminaries}\label{section 3}
\subsection{Main assumptions}
With reference to problem \eqref{1.1} we suppose that \label{section3.1}
\renewcommand{\labelenumi}{{(A\arabic{enumi})}}
\begin{enumerate}
\item $P$ and $Q$ are Carath\'eodory functions, respectively in
$\Omega\times\R$ and $\Gamma_1\times\R$, that $P(x,v)v\ge 0$ a.e. in $\Omega\times\R$, $Q(x,v)v\ge 0$ a.e. on $\Gamma_1\times\R$,
 and  there are $\alpha\in
L^\infty(\Omega)$, $\beta\in L^\infty(\Gamma_1)$, $\alpha,\beta\ge
0$, and $1<\widetilde{m}\le m$, $1<\widetilde{\mu}\le\mu$, $c_m,c_\mu\ge 0$,  such that
\begin{equation}\label{2.19}
\begin{aligned}
&|P(x,v)|\le c_m\alpha(x) (|v|^{\widetilde{m}-1}+|v|^{m-1})       \,\, \text{for a.a. $x\in\Omega$, all $v\in\R$;}\\
&|Q(x,v)|\le c_\mu\beta(x) (|v|^{\widetilde{\mu}-1}+|v|^{\mu-1}) \quad \text{for a.a. $x\in\Gamma_1$, all $v\in\R$;}
\end{aligned}
\end{equation}
\item $f$ and $g$ are Carath\'eodory functions, respectively in
$\Omega\times\R$ and  $\Gamma_1\times\R$, and there $p,q\ge 2$, $c_p,c_q\ge 0$  such that
\begin{equation}\label{2.20}
\begin{aligned}
&|f(x,u)|\le c_p(1+|u|^{p-1}),\qquad\text{for a.a. $x\in\Omega$ and all $u\in\R$,}\\
&|g(x,u)|\le c_q(1+|u|^{q-1}),\qquad\text{for a.a. $x\in\Gamma_1$ and all $u\in\R$;}
\end{aligned}
\end{equation}
\item $p\le 1+\romega/2$ or $\essinf_\Omega \alpha>0$,  $q\le 1+\rgamma/2$ or
$\essinf_{\Gamma_1}\beta>0$, and
 \begin{equation}\label{2.21}
2\le p\le 1+\romega/ \overline{m}', \qquad 2\le q\le 1+\rgamma/
\overline{\mu}',
\end{equation}
where $\overline{m}$ and $\overline{\mu}$ are given by \eqref{1.7}.
\end{enumerate}

When $P(x,v)=\alpha(x)P_0(v)$ and  $Q(x,v)=\beta(x)Q_0(v)$ with
$\alpha\in L^\infty(\Omega)$ and $\beta\in L^\infty(\Gamma_1)$,
$\alpha,\beta\ge 0$, assumption (A1) trivially holds when $P_0,Q_0\in C(\R)$, $P_0(v)v\ge 0$, $Q_0(v)v\ge0$, and there are
$1<\widetilde{m}\le m$, $1<\widetilde{\mu}\le\mu$ such that
\begin{equation}\label{2.22}
\begin{aligned}
  Q_0(v)=O(|v|^{\widetilde{m}-1}),\qquad &P_0(v)=O(|v|^{\widetilde{\mu}-1})\qquad\text{as $v\to 0$,}\\
  Q_0(v)=O(|v|^{m-1}),\qquad &P_0(v)=O(|v|^{\mu-1})\qquad\text{as $|v|\to \infty$.}
\end{aligned}
\end{equation}
Moreover, when $f(x,u)=f_0(u)$ and $g(x,u)=g_0(u)$, assumption (A2) trivially holds when $f_0,g_0\in C(\R)$ and there are $p,q\ge 2$ such that
\begin{equation}\label{2.23}
f_0(u)=O(|u|^{p-1}),\qquad g_0(u)=O(|u|^{q-1})\qquad\text{as $|u|\to \infty$}.
\end{equation}
Consequently assumptions (A1--2) hold true for the following model nonlinearities
\begin{equation}\label{2.24}
\left\{
\begin{aligned} &P(x,v)=&P_1(v):=&\quad\alpha\left(a|v|^{\widetilde{m}-2}v+ |v|^{m-2}v\right), \\
&Q(x,v)=&Q_1(v):=&\quad\beta\left(b|v|^{\widetilde{\mu}-2}v+ |v|^{\mu-2}v\right), \\
&f(x,u)=&f_1(u):=&\quad\gamma|u|^{p-2}u+\widetilde{\gamma}|u|^{\widetilde{p}-2}u+\widetilde{\gamma}',\\
&g(x,u)=&g_1(u):=&\quad\delta|u|^{q-2}u+\widetilde{\delta}|u|^{\widetilde{q}-2}u+\widetilde{\delta}',
\end{aligned}
\right.
\end{equation}
provided
\begin{equation}\label{2.25}
\begin{gathered}
a,b,\alpha,\beta \ge 0,\qquad  \gamma, \widetilde{\gamma}, \widetilde{\gamma}',\delta, \widetilde{\delta}, \widetilde{\delta}'\in\R, \\ 1<\widetilde{m}\le m,  \qquad 1<\widetilde{\mu}\le \mu,\qquad
2\le\widetilde{p}\le p, \qquad 2\le\widetilde{q}\le q.
\end{gathered}
\end{equation}
Moreover trivially assumption (A3) hold true provided
\begin{equation}\label{2.25bis}
\left\{\begin{aligned}
 p\le
\begin{cases}
1+\romega/2&\text{if $\gamma\not=0$, $\alpha=0$,}\\
1+\romega/ \overline{m}'&\text{if $\gamma\not=0$, $\alpha>0$,}
\end{cases}\quad
q\le
\begin{cases}
1+\rgamma/2&\text{if $\delta\not=0$, $\beta=0$,}\\
1+\rgamma/ \overline{\mu}'&\text{if $\delta\not=0$, $\beta>0$,}
\end{cases}\\
 \widetilde{p}\le
\begin{cases}
1+\romega/2&\text{if $\widetilde{\gamma}\not=0$, $\alpha=0$,}\\
1+\romega/ \overline{m}'&\text{if $\widetilde{\gamma}\not=0$, $\alpha>0$,}
\end{cases}\quad
\widetilde{q}\le
\begin{cases}
1+\rgamma/2&\text{if $\widetilde{\delta}\not=0$, $\beta=0$,}\\
1+\rgamma/ \overline{\mu}'&\text{if $\widetilde{\delta}\not=0$, $\beta>0$.}
\end{cases}
\end{aligned}\right.
\end{equation}

\begin{remark}\label{remark2.1} Restricting \eqref{2.24} to the case $\widetilde{\gamma}=\widetilde{\gamma}'=\widetilde{\delta}=\widetilde{\delta}'=0$ and $\gamma,\delta\ge 0$ we trivially get the nonlinearities in problem \eqref{1.2}, and \eqref{2.25} reduces to \eqref{1.3}.
Since in this case \eqref{2.25bis} trivially reads as \eqref{1.6}, we get that assumptions (A1--3) hold true provided \eqref{1.3} and \eqref{1.6} hold.
\end{remark}

We now point out, for further reference,  that by assumption (A1) it follows that there are $c_m',c_\mu'\ge 0$ such that
\begin{equation}\label{2.26}
\begin{aligned}
&|P(x,v)|\le c'_m \left[(P(x,v)v)^{\frac 1{m'}}+(P(x,v)v)^{\frac 1{\widetilde{m}'}}\right] &&\text{for a.a. $x\in\Omega$, all $v\in\R$;}\\
&|Q(x,v)|\le c'_\mu \left[(Q(x,v)v)^{\frac 1{\mu'}}+(Q(x,v)v)^{\frac 1{\widetilde{\mu}'}}\right]&&\text{for a.a. $x\in\Gamma_1$, all $v\in\R$.}
\end{aligned}
\end{equation}
Indeed, by \eqref{2.19} trivially $P(\cdot,0)\equiv 0$  and, when $|v|\le 1$,  $|P(x,v)|\le 2c_m\alpha(x)|v|^{m-1}$, or equivalently
$|P(x,v)|^{1/m}\le [2c_m\alpha(x)]^{1/m}|v|^{1/m'}$. Since $P(x,v)v\ge 0$ we then get
\begin{multline*}
|P(x,v)|=|P(x,v)|^{1/m}|P(x,v)|^{1/m'}\le [2c_m\|\alpha\|_\infty]^{1/m}|P(x,v)|^{1/m'}|v|^{1/m'}\\
=[2c_m\|\alpha\|_\infty]^{1/m}(P(x,v)v)^{1/m'}\qquad\text{when $|v|\le 1$.}$$
\end{multline*}
Exactly the same argument show that
$$
|P(x,v)|[2c_m\|\alpha\|_\infty]^{1/\widetilde{m}}(P(x,v)v)^{1/\widetilde{m}'}\qquad\text{when $|v|\ge 1$,}$$
and hence $P$ satisfies \eqref{2.26}. The same arguments apply to $Q$ as well.

\subsection{Weak solutions}
To define weak solutions of problem \eqref{1.1} we first point out the following easy result, noticing the reader that in the sequel  we shall denote
by $\widehat{P}$, $\widehat{Q}$, $\widehat{f}$ and  $\widehat{g}$ the Nemitskii operators respectively associated to $P$, $Q$, $f$ and $g$ (see \cite[Definition~2.1, p. 15]{ambrosettiprodi} or \cite[Definition~10.57, \S 10.3.4, p.370] {RenRog}).
\begin{lemma}\label{lemma2.2}Let assumptions (A1--3) hold and $u\in L^\infty(0,T;H^1)\cap W^{1,\infty}(0,T;H^0)$, $u'\in L^{2,\rho,\theta}_{\alpha,\beta}(0,T)$ for some $0<T<\infty$. Then $\xi=\widehat{f}(u)-\widehat{P}(u_t)$ and $\eta=\widehat{g}(u_{|\Gamma})-\widehat{Q}(u_t)$ are of the form \eqref{2.15}, with $\rho=\overline{m}$ and $\theta=\overline{\mu}$.
\end{lemma}
\begin{proof} We first remark that classical results on Nemitskii operators (see \cite[Theorem~2.2, p. 16]{ambrosettiprodi})
trivially extend to abstract measure spaces. Hence, by \eqref{2.19}, the Nemitskii operator associated to $P/\alpha$, this function being $\lambda'_\alpha$ -- a.e. well defined, is continuous from $L^{\overline{m}}((0,T)\times\Omega,\lambda'_\alpha)$ to $L^{\overline{m}'}((0,T)\times\Omega,\lambda'_\alpha)$. Since, by \eqref{2.1}, \eqref{2.5} and \eqref{2.6} we have $[u_t]_\alpha\in L^{\overline{m}}((0,T)\times\Omega,\lambda'_\alpha)$, we consequently get $[P(\cdot,u_t)/\alpha]_\alpha\in L^{\overline{m}'}((0,T)\times\Omega,\lambda'_\alpha)$, i.e. $\widehat{P}(u_t)=\alpha \xi_2$, $\xi_2\in L^{\overline{m}'}((0,T)\times\Omega,\lambda'_\alpha)$.

Next, when $p\le 1+\romega/2$, by \eqref{2.20} from $u\in L^\infty(0,T;H^1(\Omega))$ and Sobolev Embedding Theorem we get
$\widehat{f}(u)\in L^\infty(0,T; L^2(\Omega))$. When $p>1+\romega/2$, by assumption (A3) we have $\essinf_\Omega \alpha>0$ so
$L^{\overline{m}'}((0,T)\times\Omega,\lambda'_\alpha)=L^{\overline{m}'}((0,T)\times\Omega)$, the norms being equivalent.
By \eqref{2.20}--\eqref{2.21} we thus get $\widehat{f}(u)\in L^\infty(0,T;L^{\overline{m}'}(\Omega))$. Consequently, as $1/\alpha\in L^\infty(\Omega)$,
$\widehat{f}(u)=\alpha\xi'_2$, $\xi'_2\in L^{\overline{m}'}((0,T)\times\Omega,\lambda'_\alpha)$.
Hence in both cases $\xi$ is in the form prescribed by \eqref{2.15} with $\rho=\overline{m}$.

The same arguments show that $\eta$ is in the form prescribed by \eqref{2.15} with $\theta=\overline{\mu}$.
\end{proof}
Thanks to Lemmas~\ref{lemma2.1} and ~\ref{lemma2.2} the following definition makes sense.
\begin{definition}\label{definition2.2} Let assumptions (A1--3) hold and $U_0=(u_0,u_1)\in\cal{H}$. By a weak solution of
problem \eqref{1.1} in $[0,T]$, $0<T<\infty$,  we mean  a weak solution of
\eqref{2.14} with
\begin{equation}\label{2.27}
\xi=\widehat{f}(u)-\widehat{P}(u_t),\quad
\eta=\widehat{g}(u_{|\Gamma})-\widehat{Q}(u_t), \quad
\rho=\overline{m}\quad\text{and}\quad \theta=\overline{\mu},
\end{equation} such that $U(0)=U_0$.
By a weak solution of \eqref{1.1} in $[0,T)$, $0<T\le\infty$, we mean $u\in
L^\infty_{\text{loc}}([0,T);H^1)$  which is a weak solution of
\eqref{1.1} in $[0,T']$ for any $T'\in (0,T)$. Such a solution is
called \emph{maximal} if it has no proper extensions and \emph{global} if $T=\infty$.
\end{definition}
We now introduce the primitives of $f$ and $g$
\begin{equation}\label{2.28}
F(x,u)=\int_0^u f(x,\tau)\,d\tau,\qquad\text{and}\quad G(x,u)=\int_0^u g(x,\tau)\,d\tau,
\end{equation}
the potential functional $J: H^{1,p,q}\to\R$ given by
\begin{equation}\label{2.29}
J(u)=\int_\Omega F(\cdot,u)\,dx+\int_{\Gamma_1}G(\cdot,u)\,d\sigma,
\end{equation}
and the energy functional $\cal{E}:\cal{H}^{p,q}\to\R$ given by
\begin{equation}\label{2.30}
 \cal{E}(u,v)=\frac 12 \|v\|_{H^0}^2+\frac 12 \int_\Omega |\nabla u|^2 dx+\frac 12 \int_{\Gamma_1}|\nabla_\Gamma u|_\Gamma^2 d\sigma-J(u).
\end{equation}
By standard results on Nemitskii and potential operators (see \cite[pp. 16--22]{ambrosettiprodi})
we have $J\in C^1(H^{1,p,q})$, with Fr\`echet derivative $J'=(\widehat{f}, \widehat{g})$, and $\cal{E}\in C^1(\cal{H}^{p,q})$.

Weak solutions of \eqref{1.1} enjoy good properties, as shown in the next result.
\begin{lemma}\label{lemma2.3}
Let assumptions (A1--3) hold and $u$ be a weak solution of \eqref{1.1} in $[0,T)$. Then
\renewcommand{\labelenumi}{{\roman{enumi})}}
\begin{enumerate}
\item  $U\in C([0,T);\cal{H})$ and  the energy identity
\begin{multline}\label{2.31}
\frac 12\|u'\|_{H^0}^2+\frac 12
\int_\Omega |\nabla u|^2 dx+\frac 12 \int_{\Gamma_1}|\nabla u|_\Gamma^2 d\sigma\Bigg|_s^t+\int_s^t\left(\int_\Omega \widehat{P}(u_t)
u_t \,dx \right.\\
\left.+\int_{\Gamma_1}\widehat{Q}(u_t) u_t \,d\sigma-\int_\Omega \widehat{f}(u)
u_t \,dx-\int_{\Gamma_1}\widehat{g}(u) u_t \,d\sigma\right)d\tau=0
\end{multline}
holds for all $0\le s\le t<T$;
\item the generalized  distribution identity
\begin{multline}\label{2.32}
(u',\phi)_{H^0}\Big|_0^t
+\int_0^t\left[-(u',\phi')_{H^0}+\int_\Omega
\nabla u\nabla \phi\, dx+\int_{\Gamma_1} (\nabla_\Gamma
u,\nabla_\Gamma\phi)_\Gamma d\sigma\right.\\\left.+\int_\Omega
\widehat{P}(u_t)\phi \,dx +\int_{\Gamma_1}
\widehat{Q}(u_t)\phi \,d\sigma-\int_\Omega
\widehat{f}(u)\phi \,dx -\int_{\Gamma_1}  \widehat{g}(u)\phi\, d\sigma\right]=0
\end{multline}
holds for all $t\in [0,T)$ and $\phi\in C([0,T);H^1)\cap
C^1([0,T);H^0)\cap L^{2,\overline{m},\overline{\mu}}_{\alpha,\beta,\loc}([0,T))$;
\item when $U_0\in \cal{H}^{p,q}$ we  have $U\in C([0,T),\cal{H}^{p,q})$ and
\begin{gather}\label{2.33}
  J(u(t))-J(u(s))=\int_s^t \left(\int_\Omega \widehat{f}(u)u_t\,dx+\int_{\Gamma_1} \widehat{g}(u)u_t\,d\sigma\right)d\tau,\\
\label{2.34}
\cal{E}(U(t))-\cal{E}(U(s))+\int_s^t\left(\int_\Omega \widehat{P}(u_t)
u_t\,dx+\int_{\Gamma_1}\widehat{Q}(u_t) u_t\,d\sigma\right)d\tau=0
\end{gather}
for $0\le s \le t<T$.
\end{enumerate}
\end{lemma}
\begin{proof} Trivially i--ii) follow from Definition~\ref{definition2.2} and Lemma~\ref{lemma2.1}. To prove iii) we take $U_0\in \cal{H}^{p,q}$.
We first claim that $u\in C([0,T);H^1(\Omega)\cap L^p(\Omega))$. When $p\le \romega$, by Sobolev Embedding Theorem, there is nothing to prove, so let us take $p>\romega$. By \eqref{2.21} it immediately follows $\overline{m}>\romega$, so $\overline{m}^2-(\romega+1)\overline{m}+\romega>0$, that is $1+\romega/\overline{m}'<\overline{m}$ and, again by \eqref{2.21}, $p<\overline{m}$. Consequently, since by assumption (A3) we have $\essinf_\Omega \alpha>0$, $u_t\in L^{\overline{m}}_\loc([0,T);L^p(\Omega))$. Since $u_0\in L^p(\Omega)$ and $u(s)=u_0+\int_0^s u_t(\tau)\,d\tau$ in $L^2(\Omega)$ for $s\in [0,T)$, we get $u\in W^{1,\overline{m}}(0,t;L^p(\Omega))\hookrightarrow C([0,t];L^p(\Omega))$ for all $t\in [0,T)$, proving our claim. The same arguments also show that
$u_{|\Gamma}\in C([0,T);H^1(\Gamma_1)\cap L^q(\Gamma_1))$, so proving that  $U\in C([0,T),\cal{H}^{p,q})$.

To prove \eqref{2.33} we introduce the auxiliary exponents
\begin{equation}\label{2.35}
m_p=
\begin{cases}
2 &\text{if $p\le 1+\romega/2$},\\
\overline{m} &\text{if $p>\max\{\overline{m},1+\romega/2\}$},\\
p &\text{if $1+\romega/2<p\le \overline{m}$},
\end{cases}
\, \mu_q=\begin{cases}
2 &\text{if $q\le 1+\rgamma/2$},\\
\overline{\mu} &\text{if $q>\max\{\overline{\mu},1+\rgamma/2\}$},\\
q &\text{if $1+\rgamma/2<q\le \overline{\mu}$},
\end{cases}
\end{equation}
and we claim that $\widehat{f}\in C(H^1(\Omega)\cap L^p(\Omega); L^{m_p'}(\Omega))$.
 We first remark that, by standard properties of Nemitskii operators and Sobolev Embedding Theorem, $\widehat{f}\in C(L^p(\Omega); L^{p'}(\Omega))$ and, when $N\ge 3$ so $\romega<\infty$, $\widehat{f}\in C(H^1(\Omega); L^{\romega/(p-1)}(\Omega))$.
We now consider the three cases in \eqref{2.35}. When $p\le 1+\romega/2$ we have $\widehat{f}\in C(H^1(\Omega); L^2(\Omega))=C(H^1(\Omega); L^{m_p'}(\Omega))$.
When $p>\max\{m,1+\romega/2\}$ by \eqref{2.21} it follows that $\widehat{f}\in C(H^1(\Omega); L^{\overline{m}'}(\Omega))=C(H^1(\Omega); L^{m_p'}(\Omega))$. Finally, when $1+\romega/2<p\le \overline{m}$ we have $\widehat{f}\in C(L^p(\Omega); L^{p'}(\Omega))=C(L^p(\Omega); L^{m_p'}(\Omega))$, so proving our claim.
By the same arguments we get that $\widehat{g}\in C(H^1(\Gamma_1)\cap L^q(\Gamma_1); L^{\mu_q'}(\Gamma_1))$, and then
\begin{equation}\label{2.36}
J'=(\widehat{f},\widehat{g})\in C(H^{1,p,q}; L^{m_p'}(\Omega)\times L^{\mu_q'}(\Gamma_1)).
\end{equation}

Next we remark that, for any $t\in [0,T)$ we have  $u_t\in L^1(0,t;L^{m_p}(\Omega))$ and $u_t\in L^1(0,t;L^{\mu_q}(\Gamma_1))$ in all three cases considered in \eqref{2.35}, so
\begin{equation}\label{2.37}
u\in W^{1,1}(0,t; L^{m_p}(\Omega)\times L^{\mu_q}(\Gamma_1)).
\end{equation}
By \eqref{2.36}--\eqref{2.37} we can then apply the abstract version of the classical chain rule proved in \cite[Lemma~7.1]{Dresda2}, by taking $X_1=H^{1,p,q}$ and $Y_1=L^{m_p}(\Omega)\times L^{\mu_q}(\Gamma_1)$, from which we get that $J\cdot u\in W^{1,1}(0,t)$ and
$(J\cdot u)'=\int_\Omega \widehat{f}(u)u_t+\int_{\Gamma_1}\widehat{g}(u)u_t$ a.e. in $(0,t)$. Being $t\in [0,T)$ arbitrary \eqref{2.33} follows. Finally \eqref{2.34} simply follows by recalling \eqref{2.30} and combining \eqref{2.31} with \eqref{2.34}.
\end{proof}
\subsection{An elementary result} We conclude this section by pointing out  the following elementary result, which should be well--known, but for which
we do not have a precise reference. We also sketch its proof for the reader's convenience.
 \begin{lemma}\label{lemma3.1}Let $l>1$, $c>0$, $0<T\le\infty$ and $\psi\in W^{1,1}_\loc([0,T))$ be such that
\begin{equation}\label{A1}
\begin{cases} \psi'\ge |\psi|^l-c\quad\text{a.e. in $(0,T)$}\\
\psi(0)=\psi_0>c^{1/l}.
\end{cases}
 \end{equation}
 Then $T\le T_m(\psi_0):=\int_{\psi_0}^\infty \frac {d\tau}{\tau^l-c}<\infty$ and $\psi(t)\to\infty$ as $t\to T_m(\psi_0)^-$ provided $T=T_m(\psi_0)$.
  \end{lemma}
\begin{proof}
We first consider the Cauchy problem $y'= |y|^l-c$, $y(0)=\psi_0$. Since $\psi_0>c^{1/l}$, by standard ODE's Theory and separation of variables, it has a unique maximal classical solution $y\in C^2(-\infty,T_m(\psi_0))$  given by $y(t)=B^{-1}(t)$,  where $B:(c^{1/l},\infty)\to (-\infty,T_m(\psi_0))$ is strictly increasing and surjective, so $y(t)\to\infty$ as  $t\to T_m(\psi_0)^-$. Then, since the standard comparison argument for ODE's (see for example \cite[Chapter 1, Theorem 1.3, p. 27]{teschl}) trivially extends to generalized solutions, since the function $y\mapsto |y|^l-c$ is locally Lipschitz continuous, and since $y'-|y|^l-c\le \psi'-|\psi|^l-c$ and $y(0)=\psi(0)$ by \eqref{A1}, by comparison we get $y\le \psi$ in $[0,T)$, from which the entire statement follows.
\end{proof}
\section{Global nonexistence with two sources}\label{section4}
In this section we state and prove our first main global nonexistence result for weak solutions of \eqref{1.1} under the following additional specific assumptions, which clearly imply that $f\not\equiv 0$ and  $g\not\equiv 0$.
\begin{itemize}
\item[(F1)] There are $\gamma_0>0$ and $\gamma_1\ge 0$ such that
$$f(x,u)u-2F(x,u)\ge \gamma_0\,|u|^p-\gamma_1\qquad\text{for a.a. $x\in\Omega$ and all $u\in\R$;}$$
\item[(G1)] there are $\delta_0>0$ and $\delta_1\ge 0$ such that
$$g(x,u)u-2G(x,u)\ge \delta_0\,|u|^q-\delta_1\qquad\text{for a.a. $x\in\Gamma_1$ and all $u\in\R$.}$$
\end{itemize}
We now check that, when (A2) holds,  assumptions (F1) and (G1) are respectively equivalent to the following ones:
\begin{itemize}
\item[(F1)$'$] there are $\gamma_2>0$ and $M_f\ge 0$ such that
$$f(x,u)u-2F(x,u)\ge \gamma_2\,|u|^p\qquad\text{for a.a. $x\in\Omega$ and all $|u|\ge M_f$;}$$
\item[(G1)$'$] there are $\delta_2>0$ and $M_g\ge 0$ such that
$$g(x,u)u-2G(x,u)\ge \delta_2\,|u|^q\qquad\text{for a.a. $x\in\Gamma_1$ and all $|u|\ge M_g$.}$$
\end{itemize}
Indeed, when (F1) holds, we get (F1)$'$ by choosing $\gamma_2=\gamma_0/2$ and $M_f= (2\gamma_1/\gamma_0)^{1/p}$. Conversely, when (F1)$'$ holds, then (F1) is trivial when $M_f=0$. When $M_f>0$ it also follows by (F1)$'$ by  choosing $\gamma_0=\min\{\gamma_2,c_f/2M_f^p\}$ and $\gamma_1=3c_f/2$, where $c_f=3c_p(M_f+M_f^p)$. Indeed when $|u|\ge M_f$ then the inequality in (F1) holds. When  $|u|\le M_f$, by \eqref{2.20} and \eqref{2.28} one easily gets that
$|f(x,u)u-2F(x,u)|\le c_f$ for a.a. $x\in\Omega$. Since we also have $\gamma_0|u|^p-\gamma_1\le \frac 12 \gamma_0M_f^p-\frac 32 c_f\le -c_f$ we get (F1).
The equivalence between (G1) and (G1)$'$ is checked by the same arguments.

Hence, when $f(x,u)=f_0(u)$ and $g(x,u)=g_0(u)$, with $f_0,g_0\in C(\R)$ verifying \eqref{2.33}, denoting by $F_0$ and $G_0$ their primitives still defined by \eqref{2.28}, assumptions (F1) and (G1) respectively reduce to
\begin{equation}\label{3.1}
\varliminf_{|u|\to\infty} \frac {f_0(u)u-2F_0(u)}{|u|^p}>0,\quad \text{and}\quad
\varliminf_{|u|\to\infty} \frac {g_0(u)u-2G_0(u)}{|u|^q}>0.
\end{equation}
\begin{remark}\label{remark3.1}
When dealing with the model nonlinearities $f_1$ and $g_1$ defined in \eqref{2.24}, conditions \eqref{3.1} respectively hold when
\begin{equation}\label{3.2}
 \gamma>0, \quad 2\le \widetilde{p}<p, \qquad\text{and}\quad  \delta>0, \quad 2\le \widetilde{q}<q.
\end{equation}
When restricting  to the case $\widetilde{\gamma}=\widetilde{\gamma}'=\widetilde{\delta}=\widetilde{\delta}'=0$ and $\gamma,\delta\ge 0$, as in problem \eqref{1.2}, see Remark~\ref{remark2.1}, \eqref{3.2} respectively reduce to
\begin{equation}\label{3.3}
\gamma>0, \quad p>2, \qquad\text{and}\quad  \delta>0,\quad q>2,
\end{equation}
so that  assumptions (A1--3), (F1) and (G1) hold true for problem \eqref{1.2} provided \eqref{1.3}, \eqref{1.6} and \eqref{3.3} hold.
\end{remark}
Our first main global nonexistence result is the following one.
\begin{theorem}\label{theorem3.1} Let assumptions (A1--3), (F1), (G1) hold, and
\begin{equation}\label{3.4}
  p>\overline{m}, \qquad q>\overline{\mu}.
\end{equation}
Then, for any $U_0\in\cal{H}$ such that $\cal{E}(U_0)<0$ problem \eqref{1.1} does not admit global weak solutions.
\end{theorem}
\begin{proof} We first recall that, as explained in Remark~\ref{remark2.1}, by \eqref{2.21} and \eqref{3.4} it follows that $\overline{m},p<\romega$ and $\overline{\mu},q<\rgamma$, so $\cal{H}^{p,q}=\cal{H}$ in Lemma~\ref{lemma2.3}--iii).

The proof is based on a contradiction argument, so let $u$ be a global weak solution of \eqref{1.1} with $\cal{E}(U_0)<0$. We introduce the auxiliary function
\begin{equation}\label{3.4bis}
  \cal{K}(t)=-\cal{E}(U(t)).
\end{equation}
By Lemma~\ref{lemma2.3}--iii) and assumption (A1) the function $\cal{K}$ belongs to  $W^{1,1}_\loc([0,\infty))$ and
\begin{equation}\label{3.5}
\cal{K}'(t)=\int_\Omega \widehat{P}(u_t(t))u_t(t)\,dx+\int_{\Gamma_1} \widehat{Q}(u_t(t))u_t(t)\,d\sigma\ge 0\qquad \text{for a.a. $t>0$.}
\end{equation}
Consequently, since $\cal{E}(U_0)<0$, by \eqref{3.4bis} and \eqref{2.30} we have
\begin{equation}\label{3.6}
0<\cal{K}_0:=\cal{K}(0)\le \cal{K}(t)\le J(u(t))\quad \text{for all $t\ge 0$.}
\end{equation}
We now remark that, by assumption (A2), one easily gets the existence of $c_p',c_q'>0$ such that
\begin{equation}\label{3.7}
\begin{aligned}
&|F(x,u)|\le c_p'(1+|u|^p) \qquad&&\text{for a.a. $x\in\Omega$ and all $u\in\R$,}\\
&|G(x,u)|\le c_q'(1+|u|^q) \qquad&&\text{for a.a. $x\in\Gamma_1$ and all $u\in\R$.}
\end{aligned}
\end{equation}
By \eqref{2.29}, \eqref{3.5} and \eqref{3.6} we thus get
\begin{equation}\label{3.8}
0<\cal{K}_0\le \cal{K}(t)\le J_1(u(t)) \quad\text{for all $t\ge 0$,}
\end{equation}
where $J_1$ is defined by
\begin{equation}\label{3.9}
J_1(u)=c_p'|\Omega|+c_q'\sigma(\Gamma_1)+c_p'\|u\|_p^p+c_q'\|u\|_{q,\Gamma_1}^q,\qquad u\in H^1.
\end{equation}
By Lemma~\ref{lemma2.3} we can take $\phi=u$ as a test function in the generalized distribution identity \eqref{2.32}. Omitting in the sequel, for the sake of simplicity, the explicit dependence of $u$ and $u'$ on $t$, using \eqref{2.30} and \eqref{3.4bis} we get
\begin{equation}\label{3.9bis}
\begin{aligned}
\frac d{dt}(u',u)_{H^0}=&\|u'\|_{H^0}^2-\|\nabla u\|_2^2-\|\nabla_\Gamma u\|_{2,\Gamma_1}^2\! \\
+&\int_\Omega \!\!\widehat{f}(u)u\,dx+\!\!\int_{\Gamma_1} \!\!\widehat{g}(u)u\,d\sigma
-\int_\Omega \widehat{P}(u_t)u\,dx-\int_{\Gamma_1}\widehat{Q}(u_t)u\,d\sigma\\
=&2\|u'\|_{H^0}^2+2\cal{K}(t)+\int_\Omega \left[f(\cdot,u)u-2F(\cdot,u)\right]\,dx\\
&+\int_{\Gamma_1} \left[g(\cdot,u)u-2G(\cdot,u)\right]\,d\sigma-\int_\Omega \widehat{P}(u_t)u\,dx-\int_{\Gamma_1}\widehat{Q}(u_t)u\,d\sigma
\end{aligned}
\end{equation}
for all $t\ge 0$. Using assumptions (F1) and (G1) in \eqref{3.9bis} we obtain
\begin{equation}\label{3.10}
\begin{aligned}
\frac d{dt}(u',u)_{H^0} &\ge 2\|u'\|_{H^0}^2+2\cal{K}(t)+\gamma_0\|u\|_p^p +\delta_0\|u\|_{q,\Gamma_1}^q -\gamma_1|\Omega|-\delta_1\sigma(\Gamma_1)\\
 &-\int_\Omega \widehat{P}(u_t)u\,dx-\int_{\Gamma_1}\!\!\!\!\widehat{Q}(u_t)u\,d\sigma\qquad\text{for all $t\ge 0$.}
\end{aligned}
\end{equation}
In the sequel we shall introduce several positive constants
depending on $\Omega$, $\Gamma_1$, $P$, $Q$, $f$ and $g$, on the various constants appearing in the assumptions, and \emph{on the initial data $U_0$}. Since they are fixed we shall not give further notice of this dependence and we  shall denote these
constants by $c_i$, $i\in\N$.  We shall denote positive constants depending
\emph{also} on other objects $\Upsilon_1,\ldots,\Upsilon_n$ by
$C_i=C_i(\Upsilon_1,\ldots,\Upsilon_n)$, $i\in\N$.

We now preliminarily estimate from above the last two terms in the right--hand side of \eqref{3.10}. By \eqref{2.26}, H\"{o}lder inequality and \eqref{3.5}, noticing that both integrals in it are nonnegative, we get
\begin{align*}
\cal{I}_1(t):=&\int_\Omega \widehat{P}(u_t)u\le\int_\Omega |\widehat{P}(u_t)||u|\\
\le & c'_m \left[\int_\Omega
 (\widehat{P}(u_t)u_t)^{1/m'}|u|\,dx+\int_\Omega (\widehat{P}(u_t)u_t)^{1/\widetilde{m}'}|u|\,dx\right]\\
\le &c_m'\left\{[\cal{K}'(t)]^{1/m'}\|u\|_m+[\cal{K}'(t)]^{1/\widetilde{m}'}\|u\|_{\widetilde{m}}\right\}\qquad\text{for all $t\ge 0$.}
\end{align*}
Since $\widetilde{m}\le m<p$, using \eqref{3.8} and \eqref{3.9}, previous estimate yields that
\begin{equation}\label{3.10bis}
\begin{aligned}
\cal{I}_1(t)\le & c_1\left\{[\cal{K}'(t)]^{1/m'}+[\cal{K}'(t)]^{1/\widetilde{m}'}\right\}[J_1(u)]^{1/p}\\
=& c_1\left\{ [\cal{K}'(t)]^{\frac 1{m'}}[J_1(u)]^{\frac 1m}[J_1(u)]^{\frac 1p-\frac 1m}+[\cal{K}'(t)]^{\frac 1{\widetilde{m}'}}[J_1(u)]^{\frac 1{\widetilde{m}}}[J_1(u)]^{\frac 1p-\frac 1{\widetilde{m}}}\right\}\\
\le & c_2\left\{ [\cal{K}'(t)]^{\frac 1{m'}}[J_1(u)]^{\frac 1m}[\cal{K}(t)]^{\frac 1p-\frac 1m}+[\cal{K}'(t)]^{\frac 1{\widetilde{m}'}}[J_1(u)]^{\frac 1{\widetilde{m}}}[\cal{K}(t)]^{\frac 1p-\frac 1{\widetilde{m}}}\right\}
\end{aligned}
\end{equation}
for all $t\ge 0$. Setting  $k_1=1/m-1/p\in(0,1)$, and using \eqref{3.8} again, we  get
$$\cal{I}_1(t)\le c_3 \left\{ [\cal{K}'(t)]^{\frac 1{m'}}[J_1(u)]^{\frac 1m}+[\cal{K}'(t)]^{\frac 1{\widetilde{m}'}}[J_1(u)]^{\frac 1{\widetilde{m}}}\right\}[\cal{K}(t)]^{-k_1}\quad\text{for all $t\ge 0$.}$$
For any $\varepsilon\in (0,1]$, to be conveniently fixed in the sequel, using weighted Young inequality we then get
\begin{equation}\label{3.11}
\begin{aligned}
  \cal{I}_1(t)\le & c_3 \left [ (\varepsilon^{-m'}+\varepsilon^{-\widetilde{m'}})\cal{K}'(t) +(\varepsilon^{m}+\varepsilon^{\widetilde{m}})J_1(u)\right ] [\cal{K}(t)]^{-k_1}\\
  \le & c_4 \left [\varepsilon^{\widetilde{m}} J_1(u)+ \varepsilon^{-\widetilde{m'}}\cal{K}'(t) \right ] [\cal{K}(t)]^{-k_1}\quad\text{for all $t\ge 0$.}
  \end{aligned}
\end{equation}
Using exactly the same arguments and setting  $k_2=1/\mu-1/p\in(0,1)$ we get
\begin{equation}\label{3.12}
  \cal{I}_2(t):=\int_{\Gamma_1} \widehat{Q}(u_t)u\,d\sigma\le  c_5 \left [\varepsilon^{\widetilde{\mu}} J_1(u)+ \varepsilon^{-\widetilde{\mu'}}\cal{K}'(t) \right ] [\cal{K}(t)]^{-k_2}\quad\text{for all $t\ge 0$.}
\end{equation}
We combine \eqref{3.11} and \eqref{3.12} by setting $\overline{k}=\min\{k_1,k_2\}\in (0,1]$ and using \eqref{3.8} again, so getting
$$  \cal{I}_1(t)+\cal{I}_2(t)\le  c_6 \left [ (\varepsilon^{\widetilde{m}}+ \varepsilon^{\widetilde{\mu}})J_1(u) +(\varepsilon^{-\widetilde{m}'}+\varepsilon^{-\widetilde{\mu}'})\cal{K}'(t)\right ] [\cal{K}(t)]^{-\overline{k}}\quad\text{for all $t\ge 0$.}
$$
Setting $m_\sharp=\min\{\overline{m},\overline{\mu}\}>1$ the last estimate  is simplified to
\begin{equation}\label{3.13}
\cal{I}_1(t)+\cal{I}_2(t)\le  c_7 \left [ \varepsilon^{m_\sharp} J_1(u) +\varepsilon^{-m_\sharp'}\cal{K}'(t)\right] [\cal{K}(t)]^{-\overline{k}}\quad\text{for all $t\ge 0$.}
\end{equation}
After these preliminary estimates we now introduce the main Lyapunov functional
\begin{equation}\label{3.14}
  \cal{Z}(t)=[\cal{K}(t)]^{1-k}+\omega (u',u)_{H^0},\quad k\in (0,\overline{k}],\quad \omega>0,
\end{equation}
where $k$ and $\omega$ will be conveniently fixed in the sequel. By \eqref{3.10} and \eqref{3.14}
\begin{align*}
\cal{Z}'(t)=&(1-k)[\cal{K}(t)]^{-k}\cal{K}'(t)+\omega \frac d{dt}(u',u)_{H^0}\\
\ge &(1-k)[\cal{K}(t)]^{-k}\cal{K}'(t)+\omega \left [2\|u'\|_{H^0}^2+2\cal{K}(t)+\gamma_0\|u\|_p^p +\delta_0\|u\|_{q,\Gamma_1}^q\right. \\
-& \left.\gamma_1|\Omega|-\delta_1\sigma(\Gamma_1)
-\int_\Omega \widehat{P}(u_t)u\,dx-\int_{\Gamma_1}\!\!\!\!\widehat{Q}(u_t)u\,d\sigma\right]\quad\text{for all $t\ge 0$.}
\end{align*}
Using the estimate \eqref{3.13} in it and using \eqref{3.8} once again we get
\begin{align*}
\cal{Z}'(t)\ge &(1-k)[\cal{K}(t)]^{-k}\cal{K}'(t)+\omega \left [2\|u'\|_{H^0}^2+2\cal{K}(t)+\gamma_0\|u\|_p^p +\delta_0\|u\|_{q,\Gamma_1}^q\right. \\
-& \left.\gamma_1|\Omega|-\delta_1\sigma(\Gamma_1)
-c_7  \varepsilon^{m_\sharp} J_1(u)[\cal{K}(t)]^{-\overline{k}} -c_7 \varepsilon^{-m_\sharp'}\cal{K}'(t) [\cal{K}(t)]^{-\overline{k}}\right]\\
\ge &(1-k)[\cal{K}(t)]^{-k}\cal{K}'(t)+2\omega \|u'\|_{H^0}^2\!\!+2\omega\cal{K}(t)+\omega\gamma_0\|u\|_p^p +\omega\delta_0\|u\|_{q,\Gamma_1}^q \!\!-\omega \gamma_1|\Omega|\\
-&\omega\delta_1\sigma(\Gamma_1)
-c_7\omega \cal{K}_0^{k-\overline{k}} \varepsilon^{m_\sharp} J_1(u)[\cal{K}(t)]^{-k} -c_7 \omega\cal{K}_0^{k-\overline{k}}\varepsilon^{-m_\sharp'}\cal{K}'(t) [\cal{K}(t)]^{-k}
\end{align*}
for all $t\ge 0$. Hence, setting $C_1=C_1(k)=c_7 \cal{K}_0^{k-\overline{k}}$ and reordering
\begin{align*}
\cal{Z}'(t)\ge &(1-k-\omega C_1\varepsilon^{-m_\sharp'})[\cal{K}(t)]^{-k}\cal{K}'(t)+2\omega \|u'\|_{H^0}^2\!\!+2\omega\cal{K}(t)+\omega\gamma_0\|u\|_p^p \\
+&\omega\delta_0\|u\|_{q,\Gamma_1}^q -\omega \left[\gamma_1|\Omega|+\delta_1\sigma(\Gamma_1)\right]
-\omega C_1 \varepsilon^{m_\sharp} J_1(u)[\cal{K}(t)]^{-k} \quad\text{for all $t\ge 0$.}
\end{align*}
Using \eqref{3.8}--\eqref{3.9} in the last estimate we get
\begin{equation}\label{3.15}
\begin{aligned}
\cal{Z}'(t)\ge &(1-k-\omega C_1\varepsilon^{-m_\sharp'})[\cal{K}(t)]^{-k}\cal{K}'(t)+\omega \|u'\|_{H^0}^2\!\!+\omega\cal{K}(t)+\omega\cal{K}_0 \\
+&\omega\gamma_0\|u\|_p^p+\omega\delta_0\|u\|_{q,\Gamma_1}^q -\omega \left[\gamma_1|\Omega|+\delta_1\sigma(\Gamma_1)\right]
-\omega C_1 \varepsilon^{m_\sharp} \cal{K}_0^{-k}J_1(u)\\
= &(1-k-\omega C_1\varepsilon^{-m_\sharp'})[\cal{K}(t)]^{-k}\cal{K}'(t)+\omega \|u'\|_{H^0}^2\!\!+\omega\cal{K}(t)\\
+&\omega\left[\gamma_0- C_1 \varepsilon^{m_\sharp} \cal{K}_0^{-k}c_p'\right]\|u\|_p^p+\omega\left[\delta_0- C_1 \varepsilon^{m_\sharp} \cal{K}_0^{-k}c_q'\right]\|u\|_{q,\Gamma_1}^q \\
+&\omega \left\{\cal{K}_0-C_1 \varepsilon^{m_\sharp} \cal{K}_0^{-k} \left[\gamma_1|\Omega|+\delta_1\sigma(\Gamma_1)\right]\right\}-\omega c_8
 \quad\text{for all $t\ge 0$.}
\end{aligned}
\end{equation}
Remembering that $C_1=C_1(k)$ we now choose $\varepsilon=\varepsilon_0(k)\in (0,1]$ so small that
\begin{gather*}
\cal{K}_0-C_1(k) \varepsilon^{m_\sharp} \cal{K}_0^{-k} \left[\gamma_1|\Omega|+\delta_1\sigma(\Gamma_1)\right]\ge 0,\\
\gamma_0- C_1(k) \varepsilon^{m_\sharp} \cal{K}_0^{-k}c_p'\ge \gamma_0/2,\qquad \delta_0- C_1(k) \varepsilon^{m_\sharp} \cal{K}_0^{-k}c_q'\ge \delta_0/2.
\end{gather*}
With this choice \eqref{3.15} yields
\begin{align*}
\cal{Z}'(t)\ge &[1-k-\omega C_1(k)\varepsilon^{-m_\sharp'}][\cal{K}(t)]^{-k}\cal{K}'(t)+\omega \|u'\|_{H^0}^2\!\!+\omega\cal{K}(t)\\
+&\tfrac 12 \omega\gamma_0\|u\|_p^p+\tfrac 12\omega\delta_0\|u\|_{q,\Gamma_1}^q -\omega c_8
 \quad\text{for all $t\ge 0$.}
\end{align*}
We now take $\omega\in (0,\omega_1(k)]$, where $\omega_1(k):=(1-k)\, \varepsilon^{m_\sharp'}\,C_1^{-1}(k)$, so that  $1-k-\omega C_1(k)\varepsilon^{-m_\sharp'}\le 0$ and, as $\cal{K}'\ge 0$, previous estimate yields
$$
\cal{Z}'(t)\ge \omega \left[\|u'\|_{H^0}^2+\cal{K}(t)+\tfrac 12 \gamma_0\|u\|_p^p+\tfrac 12\delta_0\|u\|_{q,\Gamma_1}^q\right] -\omega c_8
 \quad\text{for all $t\ge 0$.}
$$
Hence, taking $c_9=\min\{1, \gamma_0/(2c_p'),\delta_0/(2c_q'), \cal{K}_0/2(c_p'|\Omega|+c_q'\sigma(\Gamma_1)\}>0$, we have the main estimate
\begin{equation}\label{3.16}
\cal{Z}'(t)\ge \omega c_9\left[\|u'\|_{H^0}^2+\cal{K}(t)+J_1(u)\right] -\omega c_8
 \quad\text{for all $t\ge 0$.}
\end{equation}
Since $k<1$, by \eqref{3.14} there is $\omega_2(k)>0$ such that for all $\omega\in (0,\omega_2(k)]$ we have
\begin{equation}\label{3.17}
  \cal{Z}(0)=\cal{K}_0^{1-k}+\omega (u_1,u_0)_{H^0}\ge \omega ^{1-k}c_8^{1-k}.
\end{equation}
We then choose $\omega=\omega_0(k)=\min\{\omega_1(k),\omega_2(k)\}$ and set $C_2=C_2(k)=c_9\omega_0(k)$, $C_3=C_3(k)=c_8\omega_0(k)$. In this way we can combine \eqref{3.16}--\eqref{3.17} to get
\begin{equation}\label{3.18}
\begin{cases}
\cal{Z}'(t)\ge \omega C_2\left[\|u'\|_{H^0}^2+\cal{K}(t)+J_1(u)\right] -C_3,\quad t\ge 0\\
\cal{Z}(0)\ge C_3^{1-k}
\end{cases}
\end{equation}
for all $k\in (0,\overline{k}]$.

We are now going to estimate from above $|\cal{Z}(t)|^l$, where $l=l(k)=1/(1-k)$ and $k\in (0,\overline{k}]$, so that $l\in (1,\overline{l}]$ where $\overline{l}=l(\overline{k})$. By \eqref{3.14} we have
$$|\cal{Z}(t)|^l\le \left[\cal{K}^{1-k}(t)+\omega_0\|u'\|_{H^0}\|u\|_{H^0}\right]^l\le 2^{l-1}
\left[\cal{K}(t)+\omega_0^l\|u'\|_{H^0}^l\|u\|_{H^0}^l\right].$$
Consequently, when we also take $k<1/2$, using Young inequality with conjugate exponents $2(1-k)$ and $2(1-k)/(1-2k)$, we have
\begin{equation}\label{3.19}
  |\cal{Z}(t)|^l\le 2^{l-1}
\left[\cal{K}(t)+\omega_0^2\|u'\|_{H^0}^2+\|u\|_{H^0}^{2/(1-2k)}\right].
\end{equation}
We finally choose $k=k_0:=\min\{\overline{k}, 1/2-1/p, 1/2-1/q\}\in(0,1/2)$, so that
$$\|u\|_{H^0}^{\frac 2{1-2k_0}}=\left(\|u\|_2^2+\|u\|_{2,\Gamma_1}^2\right)^{\frac 1{1-2k_0}}\le
2^{\frac {2k_0}{1-2k_0}}\left(\|u\|_2^{\frac 2{1-2k_0}}+\|u\|_{2,\Gamma_1}^{\frac 2{1-2k_0}}\right).$$
Since $2/(1-2k_0)\le p$ and $2/(1-2k_0)\le q$ and $\Omega$ is bounded the previous estimates yields
\begin{equation}\label{3.20}
\|u\|_{H^0}^{\frac 2{1-2k_0}}\le
c_{10}\left(1+\|u\|_p^p+\|u\|_{q,\Gamma_1}^q\right).
\end{equation}
Denoting $l_0=l(k_0)$, by \eqref{3.19}--\eqref{3.20} we thus have
$$
|\cal{Z}(t)|^{l_0}\le 2^{l_0-1}
\left[\cal{K}(t)+\omega_0^2\|u'\|_{H^0}^2+c_{10}\left(1+\|u\|_p^p+\|u\|_{q,\Gamma_1}^q\right)\right],
$$
and consequently, by \eqref{3.9}, we get
\begin{equation}\label{3.21}
|\cal{Z}(t)|^{l_0}\le c_{11}
\left[\cal{K}(t)+\|u'\|_{H^0}^2+J_1(u)\right].
\end{equation}
By combining the estimates \eqref{3.18} and \eqref{3.21} we get that $\cal{Z}$ satisfies the assumptions of Lemma \eqref{lemma3.1}, which gives the required contradiction.
\end{proof}
\section{Global nonexistence with interior source}\label{section5}
This section is devoted to our second main global nonexistence result for weak solutions of \eqref{1.1} when the interior source is present in it, while $g$ may also vanish. In particular we shall keep the specific assumption (F1) on $f$, while assumption (G1) is replaced by the following one.
\begin{itemize}
\item[(G2)] There is $\overline{q}>2$ such that
$$g(x,u)u\ge \overline{q}G(x,u)\ge 0\qquad\text{for a.a. $x\in\Gamma_1$ and all $u\in\R$.}$$
\end{itemize}
\begin{remark}\label{remark4.1}
When dealing with the model nonlinearity $g_1$ defined in \eqref{2.24}, assumption (G2) reduces to
$$ \delta,\widetilde{\delta}\ge 0,\qquad \widetilde{\delta}'=0,\qquad 2<\widetilde{q}\le q,$$
and then, comparing with \eqref{3.2}, assumptions (G1) and (G2) are unrelated. Moreover, by Remark \ref{remark3.1}, the couple of assumptions (F1), (G2) holds for the model nonlinearities $f_1$ and $g_1$ in \eqref{2.24} when
\begin{equation}\label{4.1}
 \gamma>0, \quad 2\le \widetilde{p}<p, \qquad\text{and}\quad \delta,\widetilde{\delta}\ge 0,\qquad \widetilde{\delta}'=0,\qquad 2<\widetilde{q}\le q.
\end{equation}
Consequently, when restricting  to the case $\widetilde{\gamma}=\widetilde{\gamma}'=\widetilde{\delta}=\widetilde{\delta}'=0$ and $\gamma,\delta\ge 0$ as in problem \eqref{1.2}, see Remark~\ref{remark2.1},  (F1) and (G2) hold provided
\begin{equation}\label{4.2}
\gamma>0, \quad p>2, \qquad\text{and}\quad  \delta\ge 0,\quad q>2,
\end{equation}
so that  assumptions (A1--3), (F1) and (G2) hold true for problem \eqref{1.2} provided \eqref{1.3}, \eqref{1.6} and \eqref{4.2} hold.
\end{remark}
Our second main global nonexistence result is the following one.
\begin{theorem}\label{theorem4.1} Let assumptions (A1--3), (F1), (G2) hold, and
\begin{equation}\label{4.3}
  p>\overline{m}, \qquad \overline{\mu}<1+p/2.
\end{equation}
Then, for any $U_0\in\cal{H}$ such that $\cal{E}(U_0)<0$ problem \eqref{1.1} does not admit global weak solutions.
\end{theorem}
\begin{proof}
We first recall that, as explained in Remark~\ref{remark1.5}, by \eqref{2.21} and \eqref{4.3} it follows that $\overline{m},p<\romega$ and $\overline{\mu},q<\rgamma$, so $\cal{H}^{p,q}=\cal{H}$ in Lemma~\ref{lemma2.3}--iii).

The proof is a variant of the one of Theorem \ref{theorem3.1}, so we shall keep all notation in it. Also in this case we use a contradiction argument, so let $u$ be a global weak solution of \eqref{1.1} with $\cal{E}(U_0)<0$. As in the quoted proof we introduce the auxiliary function
$\cal{K}$ defined in \eqref{3.4bis} and we get that $\cal{K}$ is increasing in $[0,\infty)$ and formulas \eqref{3.5}--\eqref{3.7} hold true. By \eqref{3.6}--\eqref{3.7}, also using \eqref{2.29}, we have
\begin{equation}\label{4.4}
0<\cal{K}_0\le \cal{K}(t)\le J_2(u(t)) \quad\text{for all $t\ge 0$,}
\end{equation}
where $J_2$ is defined by
\begin{equation}\label{4.5}
J_2(u)=c_p'|\Omega|+c_p'\|u\|_p^p+\int_{\Gamma_1}G(\cdot,u)\,d\sigma,\qquad u\in H^1.
\end{equation}
Also in this case we get the identity \eqref{3.9bis} and, for any $\varepsilon\in \left(0,\frac{\overline{q}-2}2\right)$ to be fixed in the sequel, using \eqref{2.30} and \eqref{3.4bis}, we rewrite it as
\begin{align*}
\frac d{dt}(u',u)_{H^0}
=&\tfrac{4+\varepsilon}2\|u'\|_{H^0}^2+\tfrac \varepsilon 2\|\nabla u\|_2^2+\tfrac \varepsilon 2\|\nabla_\Gamma u\|_{2,\Gamma_1}^2+(2+\varepsilon)\cal{K}(t)\\
+&\int_\Omega \left[f(\cdot,u)u-(2+\varepsilon)F(\cdot,u)\right]\,dx
+\int_{\Gamma_1} \left[g(\cdot,u)u-(2+\varepsilon) G(\cdot,u)\right]\,d\sigma\\
-&\int_\Omega \widehat{P}(u_t)u\,dx-\int_{\Gamma_1}\widehat{Q}(u_t)u\,d\sigma\qquad\text{for all $t\ge 0$.}
\end{align*}
Consequently, using \eqref{3.7} and the assumptions (F1), (G2), as $\overline{q}-2-\varepsilon>\varepsilon$ we get the estimate
\begin{align*}
\frac d{dt}(u',u)_{H^0}
\ge &2\|u'\|_{H^0}^2+\frac \varepsilon 2\|\nabla u\|_2^2+2\cal{K}(t)+\gamma_0\|u\|_p^p-\gamma_1|\Omega|\\
-& \varepsilon\int_\Omega F(\cdot,u)\,dx-\int_\Omega \widehat{P}(u_t)u\,dx+\int_{\Gamma_1} \left[g(\cdot,u)u-\overline{q} G(\cdot,u)\right]\,d\sigma\\
+&(\overline{q}-2-\varepsilon)\int_{\Gamma_1} G(\cdot,u)\,d\sigma
-\int_{\Gamma_1}\widehat{Q}(u_t)u\,d\sigma\\
\ge &2\|u'\|_{H^0}^2+\frac \varepsilon 2\|\nabla u\|_2^2+2\cal{K}(t)+(\gamma_0-\varepsilon c_p')\|u\|_p^p-(\gamma_1+\varepsilon c_p')\Omega|\\
+&\varepsilon\int_{\Gamma_1} G(\cdot,u)\,d\sigma-\int_\Omega \widehat{P}(u_t)u\,dx-\int_{\Gamma_1}\widehat{Q}(u_t)u\,d\sigma.
\end{align*}
Restricting to $\varepsilon\in (0,\varepsilon_1)$, where $\varepsilon_1:=\min\left\{\frac{\overline{q}-2}2, \frac{\gamma_0}{2c_p'}\right\}$ we then get
\begin{equation}\label{4.6}
\begin{aligned}
\frac d{dt}(u',u)_{H^0}
\ge &2\|u'\|_{H^0}^2+\tfrac \varepsilon 2\|\nabla u\|_2^2+2\cal{K}(t)+\tfrac 12 \gamma_0\|u\|_p^p-(\gamma_1+\varepsilon c_p')\Omega|\\
+&\varepsilon\int_{\Gamma_1} G(\cdot,u)\,d\sigma-\int_\Omega \widehat{P}(u_t)u\,dx-\int_{\Gamma_1}\widehat{Q}(u_t)u\,d\sigma.
\end{aligned}
\end{equation}
Also in this case we estimate from below the last two terms in the right--hand side of \eqref{4.6}. Since, by \eqref{4.5} and assumption (G2) we have $c_p'\|u\|_p^p\le J_2(u)$ we can estimate the term $\int_\Omega \widehat{P}(u_t)u$ exactly as in the proof of Theorem \ref{theorem3.1}, with $J_2$ replacing $J_1$. In this way, for $\varepsilon\in (0,\varepsilon_2)$, where $\varepsilon_2=\min\{1,\varepsilon_1\}$, we get
\begin{equation}\label{4.7}
\begin{aligned}
  \cal{I}_1(t):=\int_\Omega \widehat{P}(u_t)u\,dx\le  c_4 \left [\varepsilon^{\widetilde{m}} J_2(u)+ \varepsilon^{-\widetilde{m'}}\cal{K}'(t) \right ] [\cal{K}(t)]^{-k_1}\quad\text{for all $t\ge 0$,}
  \end{aligned}
\end{equation}
where $k_1=1/m-1/p\in (0,1)$.

The estimate of the last term in the right -- hand side of \eqref{4.6} in this case is  different. Indeed, as in the estimate of $\cal{I}_1(t)$ we get
$$\cal{I}_2(t):=\int_{\Gamma_1} \widehat{Q}(u_t)u\,d\sigma\le c_\mu'\left\{[\cal{K}'(t)]^{1/\mu'}\|u\|_{\mu,\Gamma_1}+[\cal{K}'(t)]^{1/\widetilde{\mu}'}
\|u\|_{\widetilde{\mu},\Gamma_1}\right\}.$$
Since $\widetilde{\mu}\le\mu\le\overline{\mu}$ and $\sigma(\Gamma_1)<\infty$ previous estimate
yields
\begin{equation}\label{4.8}
\cal{I}_2(t)\le c_{12}\left\{[\cal{K}'(t)]^{1/\mu'}+[\cal{K}'(t)]^{1/\widetilde{\mu}'}\right\}
\|u\|_{\overline{\mu},\Gamma_1}.
\end{equation}
To estimate the term $\|u\|_{\overline{\mu},\Gamma_1}$ in \eqref{4.8} we set
$s=\frac 1{\overline{\mu}}\left(\frac 12 +\frac{p-\overline{\mu}}{p-2}\right)$. Since, by \eqref{4.3}, $p>2(\overline{\mu}-1)$, we have $(p-\overline{\mu})/(p-2)>1/2$. Hence, as $\overline{\mu}\ge 2$,
\begin{equation}\label{4.6BIS}
\frac 1{\overline{\mu}}<s<\frac {2(p-\overline{\mu})}{\overline{\mu}(p-2)}\le 1.
\end{equation}
We can then use the Trace Inequality \cite[Theorem 1.5.1.2, p. 37]{grisvard} to estimate
\begin{equation}\label{4.7BIS}
\|u\|_{\overline{\mu},\Gamma_1}\le c_{13}\|u\|_{W^{s,\overline{\mu}}(\Omega)}.
\end{equation}
Now, by interpolation theory, see \cite[Theorem 4.3.1.2 p. 317 and Remark 2, formula (9), \S 2.4.2 p. 185]{triebel}, we know that
\begin{equation}\label{4.8BIS}
  W^{s,\theta}(\Omega)=B_{\theta,\theta}^s(\Omega)= (H^1(\Omega),L^p(\Omega))_{1-s,\theta}
\end{equation}
where $(\cdot,\cdot)_{1-s,\theta}$ denotes the real interpolator functor, provided
\begin{equation}\label{4.9}
  \frac 1\theta=\frac s2+\frac {1-s}p,\qquad\text{i.e.}\quad \theta=\frac {2p}{s(p-2)+2}.
\end{equation}
Now, by \eqref{4.6} and \eqref{4.9} we have $\overline{\mu}<\theta$. Hence, as $\Omega$ is bounded and $p>2$, combining \eqref{4.7BIS} with the interpolation inequality (see \cite[Theorem 3.2.2 p. 43]{bergh}) which follows from \eqref{4.8BIS} we get
$$\|u\|_{\overline{\mu},\Gamma_1}\le c_{14}\|u\|_{H^1(\Omega)}^s\|u\|_p^{1-s}
\le c_{15}
\left(\|u\|_p+\|\nabla u\|_2^s\|u\|_p^{1-s}\right).
$$
Plugging it into \eqref{4.8} and using \eqref{4.5} and assumption (G2)  we get
\begin{equation}\label{4.11}
\begin{aligned}
\cal{I}_2(t)\le
&\quad c_{16}\left\{[\cal{K}'(t)]^{1/\mu'}+[\cal{K}'(t)]^{1/\widetilde{\mu}'}\right\}\left(\|u\|_p+\|\nabla u\|_2^s\|u\|_p^{1-s}\right)\\
\le &\quad\cal{I}_2^1(t)+\cal{I}_2^2(t),
\end{aligned}
\end{equation}
where
\begin{align}\label{4.12}
\cal{I}_2^1(t):=&c_{17}\left\{[\cal{K}'(t)]^{1/\mu'}+[\cal{K}'(t)]^{1/\widetilde{\mu}'}\right\} [J_2(u)]^{1/p},\\\label{4.12bis}
\cal{I}_2^2(t):=&c_{17}\left\{[\cal{K}'(t)]^{1/\mu'}+[\cal{K}'(t)]^{1/\widetilde{\mu}'}\right\} \|\nabla u\|_2^s \,\,[J_2(u)]^{\frac {1-s}p}.
\end{align}
We now separately estimate $\cal{I}_2^1(t)$ and $\cal{I}_2^2(t)$. We estimate the first one  exactly as we estimated the term $\cal{I}_1(t)$ in \eqref{3.10bis}, by respectively  replacing $m$, $\widetilde{m}$ and $J_1$ with $\mu$, $\widetilde{\mu}$ and $J_2$. In this way, setting $k_2=1/\mu-1/p\in (0,1)$, we get the estimate
\begin{equation}\label{4.13}
  \cal{I}_2^1(t)\le c_{18}\left [\varepsilon^{\widetilde{\mu}} J_2(u)+ \varepsilon^{-\widetilde{\mu'}}\cal{K}'(t) \right ] [\cal{K}(t)]^{-k_2}\quad\text{for all $t\ge 0$.}
\end{equation}
To estimate the term $\cal{I}_2^2(t)$ we first set
\begin{equation}\label{4.14}
k_3=\frac 1{\overline{\mu}}-\left(\frac s2+\frac{1-s}p\right)=\frac 1{\overline{\mu}}-\frac 1\theta>0.
\end{equation}
By \eqref{4.4}, \eqref{4.12} and \eqref{4.14}, since $\widetilde{\mu}\le\mu\le\overline{\mu}$,
\begin{equation}\label{4.15}
\begin{aligned}
\cal{I}_2^2(t)=&c_{17}\left\{[\cal{K}'(t)]^{\frac 1{\mu'}}\|\nabla u\|_2^s \,\,[J_2(u)]^{\frac 1\mu-\frac s2}\right.\\
+&\left.[\cal{K}'(t)]^{\frac 1{\widetilde{\mu}'}}\|\nabla u\|_2^s [J_2(u)]^{\frac 1{\overline{\mu}}-\frac s2}\right\}[J_2(u)]^{-k_3}\\
\le &c_{18}\left\{[\cal{K}'(t)]^{\frac 1{\mu'}}\|\nabla u\|_2^s \,\,[J_2(u)]^{\frac 1\mu-\frac s2}\right.\\
+&\left.[\cal{K}'(t)]^{\frac 1{\widetilde{\mu}'}}\|\nabla u\|_2^s [J_2(u)]^{\frac 1{\overline{\mu}}-\frac s2}\right\}[\cal{K}(t)]^{-k_3}\quad\text{for all $t\ge 0$.}
\end{aligned}
\end{equation}
 By \eqref{4.6BIS} we have $\frac s2<\frac {p-\overline{\mu}}{\overline{\mu}(p-2)}<\frac 1{\overline{\mu}}$, so
$$  \frac 1{\widetilde{\mu}}-\frac s2\ge \frac 1\mu-\frac s2\ge \frac 1{\overline{\mu}}-\frac s2>0.
$$
Hence, using triple Young inequality with exponents $\mu'$, $2/s$ and $1/\left(\frac 1\mu-\frac s2\right)$ and with exponents $\widetilde{\mu}'$, $2/s$ and $1/\left(\frac 1{\widetilde{\mu}}-\frac s2\right)$  we get the estimates
\begin{align}\label{4.17}
[\cal{K}'(t)]^{\frac 1{\mu'}}\|\nabla u\|_2^s [J_2(u)]^{\frac 1\mu-\frac s2}=&[\varepsilon^{-\frac {2\mu'}\mu}\cal{K}'(t)]^{\frac 1{\mu'}}\left(\varepsilon\|\nabla u\|_2\right)^s \left\{\varepsilon^2 [J_2(u)]\right\}^{\frac 1\mu-\frac s2}\\
\nonumber\le & \varepsilon^2 J_2(u)+\varepsilon^2 \|\nabla u\|_2^2+ \varepsilon^{-\frac {2\mu'}\mu}\cal{K}'(t),\\\label{4.18}
[\cal{K}'(t)]^{\frac 1{\widetilde{\mu}'}}\|\nabla u\|_2^s [J_2(u)]^{\frac 1{\widetilde{\mu}}-\frac s2}
\le & \varepsilon^2 J_2(u)+\varepsilon^2 \|\nabla u\|_2^2+ \varepsilon^{-\frac {2\widetilde{\mu}'}{\widetilde{\mu}}}\cal{K}'(t)
\end{align}
for all $t\ge 0$. Since $\varepsilon\le 1$ and $\widetilde{\mu}\le\mu$, by \eqref{4.15}, \eqref{4.17} and \eqref{4.18}
\begin{equation}\label{4.19}
\cal{I}_2^2(t)\le c_{19} \left[\varepsilon^2 J_2(u)+\varepsilon^2 \|\nabla u\|_2^2+ \varepsilon^{-\frac {2\mu'}\mu}\cal{K}'(t)\right][\cal{K}(t)]^{-k_3}\quad\text{for all $t\ge 0$.}
\end{equation}
To estimate $\cal{I}_2(t)$ we now remark that, being $p>2$, using \eqref{4.14}, we have
 $k_2=\frac 1\mu-\frac 1p\ge \frac 1{\widetilde{\mu}}-\frac 1p\ge k_3$.  By \eqref{4.4} we can combine \eqref{4.11}, \eqref{4.13} and \eqref{4.19} to get
\begin{equation}\label{4.20}
\cal{I}_2(t)\le c_{20} \left[\varepsilon^{\min\{2,\widetilde{\mu}\}} J_2(u)+\varepsilon^2 \|\nabla u\|_2^2+ \varepsilon^{-\min\left\{\widetilde{\mu}',\frac {2\mu'}\mu\right\}}\cal{K}'(t)\right][\cal{K}(t)]^{-k_3}
\end{equation}
for all $t\ge 0$. We then set $k_4=\min\{k_1,k_3\}$ and, using \eqref{4.4} again we combine \eqref{4.7} and \eqref{4.20} to get
$$\cal{I}_1(t)+\cal{I}_2(t)\le c_{21} \left[\varepsilon^{\min\{2,\widetilde{\mu},\widetilde{m}\}} J_2(u)+\varepsilon^2 \|\nabla u\|_2^2+ \varepsilon^{-\min\left\{\widetilde{\mu}',\widetilde{m}',\frac {2\mu'}\mu\right\}}\cal{K}'(t)\right][\cal{K}(t)]^{-k_4}
$$
for all $t\ge 0$. Hence, using \eqref{4.4}--\eqref{4.5}, for any $k\in (0,k_4]$ there is $C_4=C_4(k)$ such that
\begin{multline*}
  \cal{I}_1(t)+\cal{I}_2(t)\le c_{22} \varepsilon^{\min\{2,\widetilde{\mu},\widetilde{m}\}}
\left(c_p'|\Omega|+c_p'\|u\|_p^p+\int_{\Gamma_1} G(\cdot,u)\,d\sigma\right)
+c_{22}\varepsilon^2 \|\nabla u\|_2^2\\
+C_4\varepsilon^{-\min\left\{\widetilde{\mu}',\widetilde{m}',\frac {2\mu'}\mu\right\}}\cal{K}'(t)[\cal{K}(t)]^{-k}\quad\text{for all $t\ge 0$.}
\end{multline*}
Plugging this estimate into \eqref{4.6} we get
\begin{equation}\label{4.22}
\begin{aligned}
\frac d{dt}(u',u)_{H^0}
\ge &2\|u'\|_{H^0}^2+\varepsilon \left(\tfrac 12-c_{22}\varepsilon\right)\|\nabla u\|_2^2+2\cal{K}(t)\\
+&\left(\tfrac {\gamma_0}2 -c_{22}\varepsilon^{\min\{2,\widetilde{\mu},\widetilde{m}\}}c_p'\right)\|u\|_p^p\\
+&\left(\varepsilon-c_{22}\varepsilon^{\min\{2,\widetilde{\mu},\widetilde{m}\}}\right)\int_{\Gamma_1} G(\cdot,u)\,d\sigma\\-&C_4\varepsilon^{-\min\left\{\widetilde{\mu}',\widetilde{m}',\frac {2\mu'}\mu\right\}}\cal{K}'(t)[\cal{K}(t)]^{-k}\\
-&(\gamma_1+\varepsilon c_p'+c_{22}\varepsilon^{\min\{2,\widetilde{\mu},\widetilde{m}\}}c_p' )|\Omega|\quad\text{for all $t\ge 0$.}
\end{aligned}
\end{equation}
Since $\widetilde{m},\widetilde{\mu}>1$ we can
finally fix $\varepsilon=\varepsilon_0'\in (0,\varepsilon_2]$ so small that
\begin{align*}
\tfrac 12-c_{22}\varepsilon\ge \tfrac 14, \qquad
\tfrac {\gamma_0}2 -c_{22}\varepsilon^{\min\{2,\widetilde{\mu},\widetilde{m}\}}c_p'\ge \tfrac {\gamma_0}4,\qquad
\varepsilon-c_{22}\varepsilon^{\min\{2,\widetilde{\mu},\widetilde{m}\}}\ge 0
\end{align*}
With this choice, using assumption (G2), from \eqref{4.22} we get
\begin{equation}\label{4.23}
\frac d{dt}(u',u)_{H^0}
\ge 2\|u'\|_{H^0}^2+2\cal{K}(t)+\tfrac {\varepsilon_0'}4 \|\nabla u\|_2^2
+\tfrac {\gamma_0}4 \|u\|_p^p
-C_5\cal{K}'(t)[\cal{K}(t)]^{-k}-c_{23}
\end{equation}
for all $t\ge 0$, where $C_5=C_5(k)$.

We now introduce, as in the proof of Theorem~\ref{theorem3.1}, the Lyapunov functional $\cal{Z}$ given by \eqref{3.14}, with $k\in (0,k_4]$ and $\omega>0$ to be conveniently fixed in the sequel.  By \eqref{4.23} we have
\begin{align*}
\cal{Z}'(t)=&(1-k)[\cal{K}(t)]^{-k}\cal{K}'(t)+\omega \frac d{dt}(u',u)_{H^0}\ge \left(1-k-C_5\omega\right)[\cal{K}(t)]^{-k}\cal{K}'(t)\\
+&2\omega \|u'\|_{H^0}^2+2\omega\cal{K}(t)+\tfrac {\varepsilon_0'\omega}4 \|\nabla u\|_2^2+\tfrac{\gamma_0\omega}4\|u\|_p^p  -\omega c_{23}\quad\text{for all $t\ge 0$.}
\end{align*}
Hence, by taking $\omega\in (0,\omega_3]$, where $\omega_3=\omega_3(k):=(1-k)-C_5(k)^{-1}$, so that $1-k-C_5(k)\omega\le 0$, we have
\begin{equation}\label{4.24}
\begin{aligned}
\cal{Z}'(t)\ge & \omega\left[2 \|u'\|_{H^0}^2+2\cal{K}(t)+\tfrac {\varepsilon_0'}4 \|\nabla u\|_2^2+\tfrac{\gamma_0}4\|u\|_p^p  - c_{23}\right]\\
\ge & c_{24}\omega\left[\|u'\|_{H^0}^2+\cal{K}(t)+\|\nabla u\|_2^2+\|u\|_p^p  \right]- c_{23}\omega\quad\text{for all $t\ge 0$.}
\end{aligned}
\end{equation}
By \eqref{3.14}, since $k<1$, there is $\omega_4=\omega_4(k)>0$ such that for all $\omega\in (0,\omega_4(k)]$
\begin{equation}\label{4.25}
  \cal{Z}(0)=\cal{K}_0^{1-k}+\omega (u_1,u_0)_{H^0}\ge \omega ^{1-k}c_{23}^{1-k}.
\end{equation}
We then choose $\omega=\omega_0'(k)=\min\{\omega_3(k),\omega_4(k)\}$ and set $C_6=C_6(k)=c_{24}\omega_0'(k)$, $C_7=C_7(k)=c_{23}\omega_0'(k)$. In this way we can combine \eqref{4.24}--\eqref{4.25} to get
\begin{equation}\label{4.26}
\begin{cases}
\cal{Z}'(t)\ge \omega C_6\left[\|u'\|_{H^0}^2+\cal{K}(t)+\|\nabla u\|_2^2+\|u\|_p^p\right] -C_7,\quad t\ge 0\\
\cal{Z}(0)\ge C_7^{1-k},
\end{cases}
\end{equation}
for all $k\in (0,k_4]$.

Also in this case we are now going to estimate from above $|\cal{Z}(t)|^l$, where $l=l(k)=1/(1-k)$ and $k\in (0,k_4]$, so that $l\in (1,1/(1-k_4)]$. Exactly as in the proof of Theorem~\ref{theorem3.1}, also taking $k<1/2$ we get
\begin{equation}\label{4.27}
  |\cal{Z}(t)|^l\le 2^{l-1}
\left[\cal{K}(t)+{\omega_0'}^2\|u'\|_{H^0}^2+\|u\|_{H^0}^{2/(1-2k)}\right].
\end{equation}
We now set $k_5=\min\{k_4, 1/2-1/p\}$, so that $k_5<1/2$ and we restrict to $k\in (0,k_5]$.
Since $k\le 1/2-1/p$ we also have $2k/(1-2k)\le (p-2)/2$ and consequently
\begin{equation}\label{4.28}
\|u\|_{H^0}^{\frac 2{1-2k}}\le
2^{\frac {2k}{1-2k}}\left(\|u\|_2^{\frac 2{1-2k}}+\|u\|_{2,\Gamma_1}^{\frac 2{1-2k}}\right).
\end{equation}
Since $2/(1-2k)\le p$ and  $\Omega$ is bounded we have
\begin{equation}\label{4.29}
  \|u\|_2^{\frac 2{1-2k}}\le c_{25} \left(1+\|u\|_p^p\right).
\end{equation}
Moreover, by standard trace and interpolation inequalities
\begin{align*}\|u\|_{2,\Gamma_1}\le &c_{26}\|u\|_{H^{1/2}(\Gamma)}\le c_{27} \|u\|_{H^{3/4}(\Omega)}\le c_{28}\|u\|_2^{1/4}\left(\|u\|_2^2+\|\nabla u\|_2^2\right)^{3/8}\\
\le &c_{28} \left(\|u\|_2+\|u\|_2^{1/4}\|\nabla u\|_2^{3/4}\right).
\end{align*}
Consequently, since $\Omega$ is bounded, there is $C_8=C_8(k)$ such that
\begin{align*}\|u\|_{2,\Gamma_1}^{\frac 2{1-2k}}\le &c_{29} \left(\|u\|_p+\|u\|_p^{1/4}\|\nabla u\|_2^{3/4}\right)^{\frac 2{1-2k}}\le
C_8 \left(\|u\|_p^{\frac 2{1-2k}}+\|u\|_p^{\frac 1{2-4k}}\|\nabla u\|_2^{\frac 3{2(1-2k)}}\right).
\end{align*}
Since $2/(1-2k)\le p$ we have $\|u\|_p^{\frac 2{1-2k}}\le 1+\|u\|_p^p$, so by the previous estimate
\begin{equation}\label{4.30}
\|u\|_{2,\Gamma_1}^{\frac 2{1-2k}}\le
C_8 \left(1+\|u\|_p^p+\|u\|_p^{\frac 1{2-4k}}\|\nabla u\|_2^{\frac 3{2(1-2k)}}\right).
\end{equation}
We finally choose $k=k_{00}$, where $k_{00}:=\min\{k_5, (1/2-1/p)/4\}>0$. By applying Young inequality with conjugate exponents  $4(1-2k_{00})/(1-8k_{00})$ and $4(1-2k_{00})/3$
 we get $\|u\|_p^{\frac 1{2-4k_{00}}}\|\nabla u\|_2^{\frac 3{2(1-2k_{00})}}\le \|u\|_p^{\frac 2{1-8k_{00}}}+\|\nabla u\|_2^2$.
Moreover, since $k_{00}\le (1/2-1/p)/4$ we have $2/(1-8k_{00})\le p$, so $\|u\|_p^{\frac 1{2-4k_{00}}}\|\nabla u\|_2^{\frac 3{2(1-2k_{00})}}\le 1+\|\nabla u\|_2^2+\|u\|_p^p$ and consequently, by \eqref{4.30}, taking $c_{29}=2C_8(k_{00})$,
\begin{equation}\label{4.31}
\|u\|_{2,\Gamma_1}^{\frac 2{1-2k_{00}}}\le
c_{29} \left(1+\|u\|_p^p+\|\nabla u\|_2^2\right).
\end{equation}
By \eqref{4.28}--\eqref{4.29} and \eqref{4.31} we then get
$\|u\|_{H^0}^{\frac 2{1-2k_{00}}}\le
c_{30}\left(1+\|u\|_p^p+\|\nabla u\|_2^2\right)$
and consequently, by \eqref{4.27}, denoting $l_{00}= l(k_{00})$,
\begin{equation}\label{4.32}
  |\cal{Z}(t)|^{l_{00}}\le c_{31}
\left[\|u'\|_{H^0}^2+\cal{K}(t)+\|\nabla u\|_2^2+\|u\|_p^p\right].
\end{equation}
By combining the estimates \eqref{4.26} and \eqref{4.32} we get that $\cal{Z}$ satisfies the assumptions of Lemma \eqref{lemma3.1}, which gives the required contradiction.
\end{proof}
\section{Blow--up results and proofs of Theorems~\ref{theorem1.3}--\ref{theorem1.4}}\label{section6}
The aim of this section is to combine the local theory in \cite{Dresda2} with Theorems~\ref{theorem3.1}--\ref{theorem4.1} to get two blow--up
results for weak solutions of \eqref{1.1}. As explained in Remarks~\ref{remark1.1},  \ref{remark1.5} and in the proofs of Theorems~\ref{theorem3.1}--\ref{theorem4.1}, when \eqref{3.4} or \eqref{4.3} holds we necessarily have $p<\romega$ and $q<\rgamma$. Hence, for the sake of simplicity, we shall recall the local theory in \cite{Dresda2} only when neither $f$ nor $g$ is super--supercritical.
\subsection{Known results}
Beside the main assumptions (A1--3)  made in \S~\ref{section3.1}, to fit with the setting in the quoted paper we shall consider in the sequel nonlinearities satisfying also the following additional structural conditions:
\renewcommand{\labelenumi}{{(A\arabic{enumi})}}
\begin{enumerate}
\setcounter{enumi}{3}
\item $P$ (respectively $Q$) is monotone increasing in $v$ for a.a. $x\in \Omega\times\R$ ($x\in \Gamma_1$), and $P$, $Q$ are coercive, that is there are  $c''_m,c''_\mu>0$,  such that
\begin{align*}
&P(x,v)v\ge c''_m\alpha(x) |v|^m        &&\text{for a.a. $x\in\Omega$, all $v\in\R$;}\\
&Q(x,v)v\ge c''_\mu\beta(x) |v|^\mu     &&\text{for a.a. $x\in\Gamma_1$, all $v\in\R$;}
\end{align*}
\item there are constants  $c''_p,c''_q\ge 0$  such that
\begin{align*}
&|f(x,u)-f(x,v)|\le c''_p|u-v|(1+|u|^{p-2}+|v|^{p-2}),\quad\text{for a.a. $x\in\Omega$ and all $u,v\in\R$,}\\
&|g(x,u)-g(x,v)|\le c''_q|u-v|(1+|u|^{q-2}+|v|^{q-2}),\quad\text{for a.a. $x\in\Gamma_1$ and all $u,v\in\R$.}
\end{align*}
\end{enumerate}
By combining the remarks made in \S\ref{section3.1} with those in \cite[\S 2.1, pp. 4893--4894]{Dresda2} we get the following conclusions.
When $P(x,v)=\alpha(x)P_0(v)$ and $Q(x,v)=\beta(x)Q_0(v)$, with
$\alpha\in L^\infty(\Omega)$ and $\beta\in L^\infty(\Gamma_1)$, $\alpha,\beta\ge 0$, assumptions (A1) and (A4) trivially hold when $P_0,Q_0\in C(\R)$ are monotone increasing and there are $1<\widetilde{m}\le m$,  $1<\widetilde{\mu}\le \mu$ such that
\begin{align*}
&\varliminf_{v\to 0} \frac {|P_0(v)|}{|v|^{m-1}}>0,  &\varliminf_{|v|\to \infty} \frac {|P_0(v)|}{|v|^{m-1}}>0, \quad&\varliminf_{v\to 0} \frac {|Q_0(v)|}{|v|^{\mu-1}}>0, &\varliminf_{|v|\to \infty} \frac {|Q_0(v)|}{|v|^{\mu-1}}>0,
\\
&\varlimsup_{v\to 0} \frac {|P_0(v)|}{|v|^{\widetilde{m}-1}}<\infty,  &\varlimsup_{|v|\to \infty} \frac {|P_0(v)|}{|v|^{m-1}}<\infty, \quad&\varlimsup_{v\to 0} \frac {|Q_0(v)|}{|v|^{\widetilde{\mu}-1}}<\infty, &\varlimsup_{|v|\to \infty} \frac {|Q_0(v)|}{|v|^{\mu-1}}<\infty.
\end{align*}
Moreover when $f(x,u)=f_0(u)$ and $g(x,u)=g_0(u)$ assumptions (A2) and (A5) trivially hold when $f_0,g_0\in C^{0,1}_\loc(\R)$ and there are $p,q\ge 2$ such that
$$f_0'(u)=O(|u|^{p-2}),\qquad g_0'(u)=O(|u|^{q-2})\qquad\text{as $|u|\to \infty$}.$$
Finally the statement of assumption (A3) is unchanged when $P,Q,f$ and $g$ are as before.
\begin{remark}\label{remark5.1} From the previous discussion it is clear that the model nonlinearities in \eqref{2.24} satisfy assumptions (A1--5)
provided \eqref{2.25}--\eqref{2.25bis} hold. Restricting \eqref{2.24} to the case $\widetilde{\gamma}=\widetilde{\gamma}'=\widetilde{\delta}=\widetilde{\delta}'=0$ and $\gamma,\delta\ge 0$ we trivially get, see also Remark~\ref{remark2.1}, that the nonlinearities in problem \eqref{1.2}, satisfy assumptions (A1--5)  provided \eqref{1.3} and \eqref{1.6} hold.
\end{remark}
Since the set of assumptions (A1--5) is (slightly) more restrictive that the set of assumptions \cite[(PQ1--3), (FG1--2), (FGQP1)]{Dresda2}, the following result is a particular case of \cite[Corollary~5.2]{Dresda2}.
\begin{theorem}[\bf Existence and continuation]\label{theorem5.1} Let assumptions (A1--5)  hold, with $p\le \romega$ and $q\le \rgamma$. Then for any $U_0=(u_0,u_1)\in\cal{H}$ problem
\eqref{1.1} possesses a maximal weak solution $u\in C([0,T_{\text{max}});H^1)\cap C^1([0,T_{\text{max}});H^0)$ for some $T_{\text{max}}\in (0,\infty]$.
Moreover, denoting $U=(u,u')$, when $T_\text{max}<\infty$ we have $\varlimsup_{t\to T_{\text{max}}^-}\|U(t)\|_{\cal{H}}= \infty$.
\end{theorem}
\begin{remark}\label{remark5.2} In \cite[Theorem 5.1 and Corollary 5.2]{Dresda2} the conclusion that when $T_\text{max}<\infty$ we have $\varlimsup_{t\to T_{\text{max}}^-}\|U(t)\|_{\cal{H}}= \infty$  was stated only for the weak maximal solution built there. On the other hand the proof of \cite[Theorem 5.1]{Dresda2} makes evident that the same conclusion holds true \emph{for any} weak maximal solution of \eqref{1.1}.
\end{remark}
Although Theorem~\ref{theorem5.1} gives the necessary motivation for studying the behavior of weak solutions of \eqref{1.1} and includes a basic continuation result, a more precise behavior as $t\to T_\text{max}^-$ (when $T_\text{max}<\infty$) is known when also \cite[assumption (FG2)$'$, p. 4898]{Dresda2} holds true.
We recall it here for the reader convenience.
\renewcommand{\labelenumi}{{(A\arabic{enumi})}}
\begin{enumerate}
\setcounter{enumi}{5}
\item
If $p>1+\romega/2$ then $N\le 4$,  $f(x,\cdot)\in C^2(\R)$ for a.a. $x\in\Omega$
and there is $c_p'''\ge 0$ such that
\begin{equation}\label{5.3}
|f_u(x,u)-f_u(x,v)|\le c'''_p |u-v|(1+|u|^{p-3}+|v|^{p-3})
\end{equation}
for a.a. $x\in\Omega$ and all $u,v\in\R$.
\\
If $q>1+\rgamma/2$ then $N\le 5$,  $g(x,\cdot)\in C^2(\R)$ for a.a. $x\in\Gamma_1$
and there is $c_q'''\ge 0$ such that
\begin{equation}\label{5.4}
|g_u(x,u)-g_u(x,v)|\le c'''_q |u-v|(1+|u|^{q-3}+|v|^{q-3})
\end{equation}
for a.a. $x\in\Gamma_1$ and all $u,v\in\R$.
\end{enumerate}
\begin{remark}It is worth noting that, when $N\le 4$, we have $1+\romega/2\ge 3$, so $p>3$ in \eqref{5.3}. Similarly $q>3$ in \eqref{5.4}.
\end{remark}
As showed in \cite[Remark~5.3, p. 4898]{Dresda2}, when $f(x,u)=f_0(u)$ and $g(x,u)=g_0(u)$ assumption (A6) reduces to the following one:

if  $p>1+\romega/2$ then $N\le 4$,
$f_0\in C^2(\R)$ and $f_0''(u)=O(|u|^{p-3})$ as $|u|\to\infty$;

if $q>1+\rgamma/2$ then $N\le 5$, $g_0\in C^2(\R)$  and $g_0''(u)=O(|u|^{q-3})$ as $|u|\to\infty$.
\begin{remark}\label{remark5.4}
From the previous discussion is then clear that the model nonlinearities in \eqref{2.24} satisfy  assumption (A6) provided

if  $\gamma\not=0$ and $p>1+\romega/2$ then $N\le 4$ and $\widetilde{p}>3$ or $\widetilde{\gamma}=0$;

if  $\delta\not=0$ and $q>1+\rgamma/2$ then $N\le 5$ and $\widetilde{q}>3$ or $\widetilde{\delta}=0$.

In particular, when $\widetilde{\gamma}=\widetilde{\gamma}'=\widetilde{\delta}=\widetilde{\delta}'=0$ and $\gamma,\delta\ge 0$ as in problem \eqref{1.2}, assumption (A6) holds when \eqref{1.14} does.
\end{remark}
Since the set of assumptions (A1--6) is (slightly) more restrictive that the set of assumptions \cite[(PQ1--3), (FG1), (FG2)$'$, (FGQP1)]{Dresda2}, the following result is a particular case of \cite[Theorem~6.2]{Dresda2}.
\begin{theorem}[\bf Existence, uniqueness, continuation]\label{theorem5.2} Let assumptions (A1--6)  hold, with $p\le \romega$ and $q\le \rgamma$. Then for any $U_0=(u_0,u_1)\in\cal{H}$ the maximal weak solution $u$ of problem \eqref{1.1}  in Theorem~\ref{theorem5.1} is unique.
Moreover, when $T_\text{max}<\infty$, we have $\lim_{t\to T_{\text{max}}^-}\|U(t)\|_{\cal{H}}= \infty$.
\end{theorem}
\subsection{Blow--up results and proofs of Theorems~\ref{theorem1.3}--\ref{theorem1.4}}
Our first blow--up result is an application of Theorem~\ref{theorem3.1}.
\begin{theorem}[\bf Blow--up with two sources]\label{theorem5.3} Let assumptions (A1--5), (F1), (G1) and \eqref{3.4} hold. Then for any $U_0=(u_0,u_1)\in\cal{H}$ such that $\cal{E}(U_0)<0$ and any  maximal weak solution  of problem \eqref{1.1}  one has
 $T_\text{max}<\infty$ and \eqref{1.26} holds. Moreover when also assumption (A6) holds we can replace $\varlimsup_{t\to T_{\text{max}}^-}$ with $\lim_{t\to T_{\text{max}}^-}$ in \eqref{1.26}.
\end{theorem}
\begin{proof} By Theorem~\ref{theorem5.1} and Remark~\ref{remark5.2} for any maximal weak solution $u$ of \eqref{1.1} either $T_\text{max}=\infty$ or
$T_\text{max}<\infty$ and $\varlimsup_{t\to T_{\text{max}}^-}\|U(t)\|_{\cal{H}}= \infty$. Theorem~\ref{theorem3.1} allows to exclude the first alternative.  To conclude the proof of \eqref{1.26} we now remark that, since $\cal{E}(U_0)<0$, by the energy identity \eqref{2.31} and assumption (A1) we have $\cal{E}(U(t))<0$ for all $t\in [0,T_\text{max})$, that is
\begin{equation}\label{5.3bis}
  \tfrac 12 \|u'(t)\|_{H^0}^2+\tfrac 12 \int_\Omega |\nabla u(t)|^2\,dx+\tfrac 12 \int_{\Gamma_1} |\nabla_\Gamma u(t)|_\Gamma^2\,d\sigma<J(u(t)).
\end{equation}
Hence, by \eqref{2.8}, \eqref{2.29} and \eqref{3.7}, since $p,q\ge 2$, we have
\begin{align*}
\tfrac 12 \|U(t)\|_{\cal{H}}^2<&\tfrac 12 \|u(t)\|_{2,\Gamma_1}^2+c_p'|\Omega|+c_q'\sigma(\Gamma_1)+c_p'\|u(t)\|_p^p+c_q'\|u(t)\|_{q,\Gamma_1}^q\\
\le & c_{32}+c_{33}\left(\|u(t)\|_p^p+\|u(t)\|_{q,\Gamma_1}^q\right),
\end{align*}
concluding the proof of \eqref{1.26}.
\end{proof}
\begin{remark}\label{remark5.5} It is useful to combine Remarks~\ref{remark3.1} and \ref{remark5.1} to point out when the model nonlinearities in \eqref{2.24} satisfy the full set of assumptions of Theorem~\ref{theorem5.3}, but for (A6). Assumptions (A1--5), (F1) and  (G1) are satisfied when \eqref{2.25} holds, with $\gamma,\delta>0$,  $\widetilde{p}<p$, $\widetilde{q}<q$,
\begin{equation}\label{1.6bis}
\begin{aligned}
 &p\le
\begin{cases}
1+\romega/2&\text{if $\alpha=0$,}\\
1+\romega/ \overline{m}'&\text{if  $\alpha>0$,}
\end{cases}
&&q\le
\begin{cases}
1+\rgamma/2&\text{if $\beta=0$,}\\
1+\rgamma/ \overline{\mu}'&\text{if $\beta>0$,}
\end{cases}\\
&\widetilde{\gamma}=0 \text{ or } \widetilde{p}\le
\begin{cases}
1+\romega/2&\text{if  $\alpha=0$,}\\
1+\romega/ \overline{m}'&\text{if $\alpha>0$,}
\end{cases}
&&\widetilde{\delta}=0 \text{ or }\widetilde{q}\le
\begin{cases}
1+\rgamma/2&\text{if  $\beta=0$,}\\
1+\rgamma/ \overline{\mu}'&\text{if  $\beta>0$.}
\end{cases}
\end{aligned}
\end{equation}
Moreover, in assumption \eqref{3.4} we can take $\overline{m}=2$ when $\alpha=0$ and $\overline{\mu}=2$ when $\beta=0$.
Hence, for the model nonlinearities in \eqref{2.24}, since $\gamma,\delta>0$, assumption \eqref{3.4} can be rewritten as
\begin{equation}\label{5.5}
p>
\begin{cases}
2&\text{if $\alpha=0$,}\\
\overline{m}&\text{if $\alpha>0$,}
\end{cases}\qquad
q>
\begin{cases}
2&\text{if $\beta=0$,}\\
\overline{m}&\text{if $\beta>0$.}
\end{cases}
\end{equation}
In particular, when $\widetilde{\gamma}=\widetilde{\gamma}'=\widetilde{\delta}=\widetilde{\delta}'=0$ and $\gamma,\delta\ge 0$ as in problem \eqref{1.2}, the full set of assumptions of Theorem~\ref{theorem5.3}, but for (A6), is satisfied when $\gamma,\delta\ge 0$, \eqref{1.3}, \eqref{1.6} and \eqref{5.5} hold true.
Finally, assumption (A6) can be checked as in Remark~\ref{remark5.4}, and in particular in problem \eqref{1.2} is satisfied when
$N\le 4$, or $N=5$ and $p\le 1+\romega/2=8/3$, or $N\ge 6$, $p\le 1+\romega/2$, $q\le 1+\rgamma/2$.
\end{remark}
Our second blow--up result is an application of Theorem~\ref{theorem4.1}.
\begin{theorem}[\bf Blow--up with interior source]\label{theorem5.4} Let assumptions (A1--5), (F1), (G2) and \eqref{4.3} hold. Then for any $U_0=(u_0,u_1)\in\cal{H}$ such that $\cal{E}(U_0)<0$ and any  maximal weak solution  of problem \eqref{1.1}  one has
 $T_\text{max}<\infty$ and
 \begin{equation}\label{5.6}
   \varlimsup_{t\to T_{\text{max}}^-}\|U(t)\|_{\cal{H}}=  \varlimsup_{t\to T_{\text{max}}^-}\|u(t)\|_p+\|u(t)\|_{2,\Gamma_1}+\int_{\Gamma_1}G(\cdot, u(t))\,d\sigma=\infty.
 \end{equation}
 Moreover when also assumption (A6) holds we can replace $\varlimsup_{t\to T_{\text{max}}^-}$ with $\lim_{t\to T_{\text{max}}^-}$ in \eqref{5.6}.
\end{theorem}
\begin{proof}By Theorem~\ref{theorem5.1} and Remark~\ref{remark5.2} for any maximal weak solution $u$ of \eqref{1.1} either $T_\text{max}=\infty$ or
$T_\text{max}<\infty$ and $\varlimsup_{t\to T_{\text{max}}^-}\|U(t)\|_{\cal{H}}= \infty$, and Theorem~\ref{theorem4.1} allows to exclude the first alternative.
To conclude the proof of \eqref{5.6} we remark that \eqref{5.3bis} holds also in this case for all $t\in [0,T_\text{max})$.
Hence, by \eqref{2.8}, \eqref{2.29} and \eqref{3.7}, since $p\ge 2$, we have
\begin{align*}
\tfrac 12 \|U(t)\|_{\cal{H}}^2<&\tfrac 12 \|u(t)\|_{2,\Gamma_1}^2+c_p'|\Omega|+c_p'\|u(t)\|_p^p+\int_{\Gamma_1} G(\cdot, u(t))\,d\sigma\\
\le & c_{34}+c_{35}\left(\|u(t)\|_p+\|u(t)\|_{2,\Gamma_1}+\int_{\Gamma_1} G(\cdot, u(t))\,d\sigma\right)^p,
\end{align*}
concluding the proof of \eqref{5.6}.
\end{proof}
\begin{remark}\label{remark5.6} Also in this case it is useful to combine Remarks~\ref{remark4.1} and \ref{remark5.1} to point out when the model nonlinearities in \eqref{2.24} satisfy the full set of assumptions of Theorem~\ref{theorem5.4}, but for (A6). Assumptions (A1--5), (F1) and  (G2) are satisfied when \eqref{2.25}, \eqref{3.2} and \eqref{1.6bis} hold.

Moreover, in assumption \eqref{4.3} we can take $\overline{m}=2$ when $\alpha=0$ and $\overline{\mu}=2$ when $\beta=0$.
Hence, for the model nonlinearities in \eqref{2.24},  assumption \eqref{4.3} can be rewritten as \eqref{1.21}.
In particular, when $\widetilde{\gamma}=\widetilde{\gamma}'=\widetilde{\delta}=\widetilde{\delta}'=0$ and $\gamma,\delta\ge 0$ as in problem \eqref{1.2}, the full set of assumptions of Theorem~\ref{theorem5.4}, but for (A6), is satisfied when $\gamma>0$, $\delta\ge 0$, \eqref{1.3}, \eqref{1.6} and \eqref{1.21} hold true.
Finally, assumption (A6) can be checked as in Remark~\ref{remark5.5}.
\end{remark}

We can finally prove the main results stated in \S~\ref{intro}.
\begin{proof}[Proof of Theorem~\ref{theorem1.3}] As seen in Remark~\ref{remark2.1}, problem \eqref{1.2} is a particular case of problem \eqref{1.1}, with $P,Q,f,g$ given by \eqref{2.24} and $\widetilde{\gamma}=\widetilde{\gamma}'=\widetilde{\delta}=\widetilde{\delta}'=0$,  $\gamma,\delta\ge 0$. Under the specific assumptions of Theorem~\ref{theorem1.3} we have $\gamma>0$ and $\delta=0$, so by Remark~\ref{remark5.6} the  assumptions of Theorem~\ref{theorem5.4} are satisfied. By observing that since $G\equiv 0$ in this case \eqref{5.6} and \eqref{1.23} are equivalent and  applying Theorem~\ref{theorem5.4} we then get the result but for the final statement, which follows by Remark~\ref{remark5.4} and the final statement of Theorem~\ref{theorem5.4}.
\end{proof}
\begin{proof}[Proof of Theorem~\ref{theorem1.4}] As in the previous proof problem \eqref{1.2} is a particular case of problem \eqref{1.1}, with the same $P,Q,f,g$. Under the specific assumptions of Theorem~\ref{theorem1.4} we have $\gamma,\delta>0$.
To get the statement we have to consider two different cases: $\overline{\mu}<q$ and $\overline{\mu}\ge q$.
In the first case assumption \eqref{1.24bis} coincides with \eqref{5.5}. Hence, recalling Remark~\ref{5.5}, we can apply Theorem~\ref{theorem5.3} and complete the proof. In the second case assumption \eqref{1.24bis} coincides with assumption \eqref{1.21}, so by Remark~\ref{remark5.6} we can apply Theorem~\ref{theorem5.4} and complete the proof.
\end{proof}

\def\cprime{$'$}
\providecommand{\href}[2]{#2}
\providecommand{\arxiv}[1]{\href{http://arxiv.org/abs/#1}{arXiv:#1}}
\providecommand{\url}[1]{\texttt{#1}}
\providecommand{\urlprefix}{URL }

\medskip
Received xxxx 20xx; revised xxxx 20xx.
\medskip

\end{document}